\documentclass[11pt]{article}
\makeatletter \@addtoreset{figure}{section} \makeatother
\makeatletter
\long\def\@makecaption#1#2{%
   \vskip 10\p@
   \setbox\@tempboxa\hbox{{#1}\ \ #2}%
   \ifdim \wd\@tempboxa >\hsize
       {#1}\ \ #2\par
   \else
       \hbox to\hsize{\hfil\box\@tempboxa\hfil}%
   \fi}
\makeatother \makeatletter
\newcommand
\myroman[1]{\@roman #1}
\makeatother

\newtheorem{thm}{Theorem}[section]
\newtheorem{cor}[thm]{Corollary}
\newtheorem{lem}[thm]{Lemma}
\newtheorem{rem}[thm]{Remark}

\newcommand{\qed}{{\hfill\rule{3pt}{7pt}}}
\def\pf{\noindent {\it Proof.} }

\def\qed{\hfill \rule{4pt}{7pt}}
\def\pf{\noindent {\it Proof} }
\title{\bf Linking Solutions for $p$-Laplace Equations with Nonlinear Boundary Conditions \\and Indefinite Weight}
\author{ Chungen Liu\footnote{Partially  supported by NFS
of China and 973 Program of
STM}\footnote{E-mail: liucg@nankai.edu.cn}~~~~~~~
~~~~~~~~~~
Youquan Zheng\footnote{E-mail: zhengyq@mail.nankai.edu.cn} \\
School of Mathematics and LPMC, Nankai University \\Tianjin 300071,
People's Republic of China\\}
\begin{document}
\date{}
\maketitle {\bf \noindent\large Abstract} We apply the linking
method for cones in normed spaces to $p$-Laplace equations with
various nonlinear boundary conditions. Some existence results are
obtained.


\section{Introduction and main results}
In this paper, we consider the following problems:\\
 \noindent Steklov
boundary problem
\begin{equation}\label{e1.1}
\left \{
\begin{array}
{ll} \Delta_{p}u = \varepsilon|u|^{p-2}u,
&\textrm{  in $\Omega$},\\
|\nabla u|^{p-2}\frac{\partial u}{\partial n}=\lambda V(x)|u|^{p-2}u
+ h(x,u),
&\textrm{  on $\partial\Omega$},\\
\end{array}
\right.
\end{equation}
\noindent No-flux boundary problem
\begin{equation}\label{e1.2}
\left \{
\begin{array}
{ll} -\Delta_{p}u + \varepsilon |u|^{p-2}u=\lambda V(x)|u|^{p-2}u +
h(x,u),
&\textrm{  in $\Omega$},\\
u=constant,
&\textrm{  on $\partial\Omega$},\\
\int_{\partial\Omega}|\nabla u|^{p-2}\frac{\partial u}{\partial
n}{\rm d}S_x=0,
\end{array}
\right.
\end{equation}
\noindent  Neumann boundary problem
\begin{equation}\label{e1.3}
\left \{
\begin{array}
{ll} -\Delta_{p}u + \varepsilon |u|^{p-2}u=\lambda V(x)|u|^{p-2}u +
h(x,u),
&\textrm{  in $\Omega$},\\
\frac{\partial u}{\partial n}=0,
&\textrm{  on $\partial\Omega$},\\
\end{array}
\right.
\end{equation}
\noindent  Robin boundary problem
\begin{equation}\label{e1.4}
\left \{
\begin{array}
{ll} -\Delta_{p}u + \varepsilon |u|^{p-2}u=\lambda V(x)|u|^{p-2}u +
h(x,u),
&\textrm{  in $\Omega$},\\
|\nabla u|^{p-2}\frac{\partial u}{\partial n}+\gamma(x)|u|^{p-2}u=0,
&\textrm{  on $\partial\Omega$}.\\
\end{array}
\right.
\end{equation}
Here $\Omega$ is a bounded domain in $\mathbf{R}^{N}$ with smooth
boundary $\partial\Omega$, ${\rm d}S_x$ is the surface element on
$\partial\Omega$, $\frac{\partial u}{\partial n}$ is the outer
normal derivative of $u$ with respect to $\partial\Omega$,
$\Delta_{p}u :=div(|\nabla u|^{p-2}\nabla u)$ is the $p$-Laplacian
operator with $p
> 1$, $\varepsilon>0$ is a constant and $V(x)\in L^{r}(\partial\Omega)$(in the case of (\ref{e1.1})) or $\in L^{r}(\Omega)$
(in the cases of (\ref{e1.2})-(\ref{e1.4})),
 where $r
= r(N,p)$ is defined by
\begin{equation}\label{e1.5}
\left \{
\begin{array}{ll}
r > (N-1)/(p-1), & if ~~1 < p < N,\\
r > 1, & if ~~p = N,\\
r = 1, & if ~~p > N.
\end{array}
\right.
\end{equation}
In problem (\ref{e1.4}), the function $\gamma(x)$ satisfies $\gamma(x)\in
L^{\infty}(\partial\Omega)$ and $\gamma(x)\ge 0$ for a.e.$x\in
\partial\Omega$.

In \cite{DL07}, the authors established and applied the linking
method for cones in normed spaces to consider the following   the
problem
\begin{eqnarray*}
\left \{
\begin{array}
{ll} -\Delta_{p}u + \varepsilon |u|^{p-2}u=\lambda V(x)|u|^{p-2}u +
h(x,u),
&\textrm{  in $\Omega$},\\
u=0,
&\textrm{  on $\partial\Omega$}\\
\end{array}
\right.
\end{eqnarray*}
 for $\lambda\in \mathbf{R}$,
$V\in L^{\infty}(\Omega)$ and $h$ satisfying ($h1'$)--($h4'$) below,
they obtained the existence of a nontrivial solution. The main goal
of this paper is to apply this method to study the problems
(\ref{e1.1})--(\ref{e1.4}).

For problem (\ref{e1.1}), we assume that $h:\partial\Omega \times
\mathbf{R} \rightarrow\mathbf{R}$ is a Carath\'{e}odory function
(i.e., $h(x,s)$ is continuous in $s$ for a.e.$x\in
\partial\Omega$ and measurable in $x$ for all $s\in \mathbf{R}$) satisfying the
following conditions:\\
(h1) if $p < N,\forall\,\epsilon >0$, $\exists\, a_{\epsilon}\in
L^{r}(\partial\Omega)$ such that $|h(x,s)| \leq
a_{\epsilon}(x)|s|^{p-1} + \epsilon|s|^{p^{*} - 1}$, $p^{*}=\frac{Np
- p}{N - p}$,
\\
$~~~~~~~$if $p = N$, $\exists \,a\in L^{r}(\partial\Omega)$, $C > 0$
and $q
> p$ such that
$|h(x,s)| \leq a(x)|s|^{p-1} + C|s|^{q-1}$, \\
$~~~~~~~$if $p
>
N$, $\forall\,S > 0$, $\exists\,a_{S}\in L^{r}(\partial\Omega)$ such
that $|h(x,s)| \leq a_{S}(x)|s|^{p-1}$
whenever $|s| \leq S$, \\
(h2) for a.e. $x\in \partial\Omega$, there hold
$\displaystyle\lim_{s\rightarrow 0}\frac{H(x,s)}{|s|^{p}} = 0$ and
$\displaystyle\lim_{|s| \rightarrow
\infty}\frac{H(x,s)}{|s|^{p}} = +\infty$, \\
(h3) there exist $\mu > p$, $\gamma_{0}\in L^{1}(\partial\Omega)$
and $\gamma_{1}\in L^{r}(\partial\Omega)$ such that $$\mu H(x,s)
\leq sh(x,s) + \gamma_{0}(x)+ \gamma_{1}(x)|s|^{p}\;{\rm ~for~
a.e.}\; x\in
\partial\Omega, {\rm and\,
every}\,s\in \mathbf{R},$$ (h4) $H(x,s) \geq 0$ for a.e.
$x\in\partial\Omega$ and every $s\in \mathbf{R}$, where $H(x,s) =
\int_{0}^{s}h(x,t){\rm d}t$.\\
For problems (\ref{e1.2})--(\ref{e1.4}), we assume $h$ satisfies the same
conditions (h1)--(h4) with $\partial \Omega$ replaced by $\Omega$
and $p^*=\frac{Np}{N-p}$.

The main results read as follow.
\begin{thm}\label{t1.1}Suppose the function $h$ satisfies the conditions (h1)--(h4)
and $V \in L^{r}(\partial\Omega)$. Then for every $\varepsilon
> 0$ and $\lambda \in \mathbf{R}$, (\ref{e1.1}) has a nontrivial
solution $u\in W^{1,p}(\Omega)$.
\end{thm}

\begin{thm}\label{t1.2}Suppose $h:\Omega\times \mathbf{R}\to
\mathbf{R}$ is a Carath\'{e}odory function satisfying the conditions
(h1)--(h4) with $\partial \Omega$ replaced by $\Omega$, $p^{*} =
\frac{Np}{N - p}$ and $V \in L^{r}(\Omega)$. Then for every
$\varepsilon
> 0$ and $\lambda \in \mathbf{R}$, the problems (\ref{e1.2}), (\ref{e1.3}) and (\ref{e1.4}) possess respectively a
nontrivial  solution $u\in W^{1,p}(\Omega)$.
\end{thm}

We note that when $\lambda \neq 0$, problems (\ref{e1.2}), (\ref{e1.3}) and
(\ref{e1.4}) are respectively equivalent to the
following problems:\\
No-flux boundary problem
\begin{eqnarray*}
\left \{
\begin{array}
{ll} -\Delta_{p}u = \lambda
(V(x)-\frac{\varepsilon}{\lambda})|u|^{p-2}u + h(x,u),
&\textrm{  in $\Omega$},\\
u=constant,
&\textrm{  on $\partial\Omega$},\\
\int_{\partial\Omega}|\nabla u|^{p-2}\frac{\partial u}{\partial
n}{\rm d}S_x=0,
\end{array}
\right.
\end{eqnarray*}
  Neumann boundary problem
\begin{eqnarray*}
\left \{
\begin{array}
{ll} -\Delta_{p}u = \lambda
(V(x)-\frac{\varepsilon}{\lambda})|u|^{p-2}u + h(x,u),
&\textrm{  in $\Omega$},\\
\frac{\partial u}{\partial n}=0,
&\textrm{  on $\partial\Omega$},\\
\end{array}
\right.
\end{eqnarray*}
  Robin boundary problem
\begin{eqnarray*}
\left \{
\begin{array}
{ll} -\Delta_{p}u  = \lambda (V(x)-
\frac{\varepsilon}{\lambda})|u|^{p-2}u + h(x,u),
&\textrm{  in $\Omega$},\\
|\nabla u|^{p-2}\frac{\partial u}{\partial n}+\gamma(x)
|u|^{p-2}u=0,
&\textrm{  on $\partial\Omega$}.\\
\end{array}
\right.
\end{eqnarray*}
Because  $V(x)- \frac{\varepsilon}{\lambda}$ is still in
$L^{r}(\Omega)$, the above three problems are exactly the following
problems
respectively.\\
  No-flux boundary problem
\begin{equation}\label{e1.6}
\left \{
\begin{array}
{ll} -\Delta_{p}u = \lambda V(x)|u|^{p-2}u + h(x,u),
&\textrm{  in $\Omega$},\\
u=constant,
&\textrm{  on $\partial\Omega$},\\
\int_{\partial\Omega}|\nabla u|^{p-2}\frac{\partial u}{\partial
n}{\rm d}S_x=0,
\end{array}
\right.
\end{equation}
 Neumann boundary problem
\begin{equation}\label{e1.7}
\left \{
\begin{array}
{ll} -\Delta_{p}u = \lambda V(x)|u|^{p-2}u + h(x,u),
&\textrm{  in $\Omega$},\\
\frac{\partial u}{\partial n}=0,
&\textrm{  on $\partial\Omega$},\\
\end{array}
\right.
\end{equation}
  Robin boundary problem
\begin{equation}\label{e1.8}
\left \{
\begin{array}
{ll} -\Delta_{p}u  = \lambda V(x)|u|^{p-2}u + h(x,u),
&\textrm{  in $\Omega$},\\
|\nabla u|^{p-2}\frac{\partial u}{\partial n}+\gamma(x)
|u|^{p-2}u=0,
&\textrm{  on $\partial\Omega$}.\\
\end{array}
\right.
\end{equation}
So Theorem 1.2 is equivalent to the following theorem (Note that the
case $\lambda = 0$ is covered by the case $\lambda \neq 0$ with $V
\equiv 0$).
\begin{thm}\label{t1.3}Suppose that the functions $h$ and $V$ satisfy the conditions as in Theorem 1.2.
Then, for every $\lambda
\in \mathbf{R}$, the problems (\ref{e1.6}), (\ref{e1.7}) and (\ref{e1.8}) possess a
nontrivial solution $u\in W^{1,p}(\Omega)$, respectively.
\end{thm}

In Theorems 1.1-1.3, if we replace  (h3) by the following condition (h5) which was introduced in \cite{J99} for $p = 2$ and in \cite{LL03} for general $p$,  the results are still true.\\
(h5) There exists a real number  $\theta \ge 1$ such that
$$\theta \mathcal{H}(x,s) \geq \mathcal{H}(x,ts) , {\rm ~for~
a.e.}\; x\in
\partial\Omega, {\rm and\,
every}\,s\in \mathbf{R}, t\in [0,1],$$
where  $\mathcal{H}(x,s):= h(x,s)s - pH(x,s)$.\\
 That is to say we have the following three results.\\
{\bf{Theorem 1.1$'$}} {\it Suppose the function $h$ satisfies the conditions (h1),(h2),(h4) and (h5), then for every $\varepsilon
> 0$ and $\lambda \in \mathbf{R}$, (\ref{e1.1}) has a nontrivial
solution $u\in W^{1,p}(\Omega)$.}\\
{\bf{Theorem 1.2$'$}} {\it Suppose $h:\Omega\times \mathbf{R}\to
\mathbf{R}$ is a Carath\'{e}odory function satisfying the conditions
(h1),(h2),(h4),(h5) with $\partial \Omega$ replaced by $\Omega$, $p^{*} =
\frac{Np}{N - p}$ and $V \in L^{r}(\Omega)$. Then for every
$\varepsilon
> 0$ and $\lambda \in \mathbf{R}$, the problems (\ref{e1.2}), (\ref{e1.3}) and (\ref{e1.4}) possess respectively a
nontrivial  solution $u\in W^{1,p}(\Omega)$.}\\
{\bf{Theorem 1.3$'$}} {\it Suppose that the functions $h$ and $V$ satisfy the conditions as in Theorem 1.2$'$.
Then, for every $\lambda
\in \mathbf{R}$, the problems (\ref{e1.6}), (\ref{e1.7}) and (\ref{e1.8}) possess a
nontrivial solution $u\in W^{1,p}(\Omega)$, respectively.}

Let $h:\overline{\Omega}\times \mathbf{R}\to \mathbf{R}$ be a
continuous function satisfying the following conditions:\\
(h$1^{'}$) if $p < N$, $\exists ~C > 0$ and $q$ satisfying $p < q <
p^{*}$, such that $|h(x,s)| \leq C(1 + |s|^{q - 1})$ ,
\\
$~~~~~~~~~~~$$p^{*}=\frac{Np}{N - p}$; if $p = N$, $\exists ~C > 0$ and $q$ satisfying $q > p$, such that $|h(x,s)| \leq C(1 + |s|^{q - 1})$; \\
$~~~~~~~~$if $p > N$, there is no restriction,\\
(h$2^{'}$) $\displaystyle\lim_{s\rightarrow 0}\frac{h(x,s)}{|s|^{p-1}} = 0$ uniformly for $x\in \overline{\Omega}$, \\
(h$3^{'}$) there exist $\mu > p$, $R > 0$ such that
$$0 < \mu H(x,s) \leq sh(x,s),\;{\rm for}\;\;|s| \geq R,$$
(h$4^{'}$) $sh(x,s) \geq 0$, where $H(x,s) =
\int_{0}^{s}h(x,t)dt$.\\ It was proved in \cite{DL07} that
(h$1^{'}$)--(h$4^{'}$) imply conditions (h1)--(h4)(with $\partial
\Omega$ replaced by $\Omega$, $p^{*} = \frac{Np}{N - p}$). So we
have the following direct consequence.
\begin{cor}\label{t1.4}Suppose $h:\overline{\Omega}\times \mathbf{R}\to
\mathbf{R}$ is a continuous function satisfying the conditions
(h$1^{'}$)--(h$4^{'}$) and $V \in L^{\infty}(\Omega)$. Then for
every $\varepsilon
> 0$ and $\lambda \in \mathbf{R}$, the problems (\ref{e1.2}), (\ref{e1.3}) and (\ref{e1.4}) possess respectively a
nontrivial  solution $u\in W^{1,p}(\Omega)$. Equivalently, for every
$\lambda \in \mathbf{R}$, the problems (\ref{e1.6}), (\ref{e1.7}) and (\ref{e1.8})
possess respectively a nontrivial  solution $u\in W^{1,p}(\Omega)$.
\end{cor}
Similarly, we have
\begin{cor}\label{t1.5}Suppose $h:\partial\Omega\times \mathbf{R}\to
\mathbf{R}$ is a continuous function satisfying the conditions
(h$1^{'}$)--(h$4^{'}$) with $\overline{\Omega}$ replaced by
$\partial\Omega$, $p^{*} = \frac{Np-p}{N - p}$ and $V \in
L^{\infty}(\partial\Omega)$. Then for every $\varepsilon
> 0$ and $\lambda \in \mathbf{R}$, (\ref{e1.1}) has a nontrivial
solution $u\in W^{1,p}(\Omega)$.
\end{cor}

The problems (\ref{e1.1})--(\ref{e1.4}), (\ref{e1.6})--(\ref{e1.8}) arise in different areas,
for example, the study of optimal constants for the Sobolev
embedding theorems(c.f.\cite{BDR04,BMR04,B06,PF01}), Non Newtonian
fluids(c.f.\cite{ADT98,AC95,AK92,D85}) and differential
geometry(c.f.\cite{E90}). Similar nonlinear boundary value problems
has been extensively studied, one can  refer to
\cite{AR09,B04,B06,BR01,C08,H04,LL,MR03,P02,P03,PBR06,W,WT06,ZZ08}
for details. In \cite{BR01}, the authors considered the following
problem
\begin{eqnarray*}
\left \{
\begin{array}
{ll} \Delta_{p}u  = |u|^{p-2}u,
&\textrm{  in $\Omega$},\\
|\nabla u|^{p-2}\frac{\partial u}{\partial n} = f(u),
&\textrm{  on $\partial\Omega$}.\\
\end{array}
\right.
\end{eqnarray*}
They proved among other cases that when $f$ has the form $\lambda
|u|^{q-2}u$ with subcritical growth, the above problem has
infinitely many solutions. In \cite{MR03}, the authors considered
the following problem
\begin{eqnarray*}
\left \{
\begin{array}
{ll} \Delta_{p}u  = |u|^{p-2}u + f(x,u),
&\textrm{  in $\Omega$},\\
|\nabla u|^{p-2}\frac{\partial u}{\partial n} = \lambda |u|^{p-2}u -
h(x,u),
&\textrm{  on $\partial\Omega$}.\\
\end{array}
\right.
\end{eqnarray*}
They obtained the existence of a solution  when $f$ and $h$ satisfy
some integral conditions of Landesmann-Laser type, and $\lambda$
equals to the first eigenvalue of the Steklov problem, i.e. the
first (minimal) $\lambda$ such that the problem
\begin{equation}\label{e1.9}\left \{
\begin{array}
{ll}\Delta_{p}u = |u|^{p-2}u,
&\textrm{  in $\Omega$},\\
|\nabla u|^{p-2}\frac{\partial u}{\partial n}=\lambda |u|^{p-2}u,
&\textrm{  on $\partial\Omega$},\\
\end{array}
\right.
\end{equation}
has a nontrivial solution. In \cite{ZZ08}, the authors considered
the following problem
\begin{eqnarray*}
\left \{
\begin{array}
{ll} -\Delta_{p}u +\lambda(x)|u|^{p-2}u = f(x,u),
&\textrm{  in $\Omega$},\\
|\nabla u|^{p-2}\frac{\partial u}{\partial n} = \eta |u|^{p-2}u,
&\textrm{  on $\partial\Omega$},\\
\end{array}
\right.
\end{eqnarray*}
where $\lambda\in L^{\infty}(\Omega)$ and $essinf_{x\in
\overline{\Omega}}\lambda(x) > 0$. They proved that if $f$ is a
superlinear and subcritical odd Carath\'{e}odory function, then the
problem they considered has infinitely many solutions for $\eta$
less than some constant. In \cite{W}, the following problem
\begin{eqnarray*}
\left \{
\begin{array}
{ll} -\Delta_{p}u  = f(x,u) - |u|^{p-2}u,
&\textrm{  in $\Omega$},\\
|\nabla u|^{p-2}\frac{\partial u}{\partial n} = \lambda |u|^{p-2}u +
g(x,u),
&\textrm{  on $\partial\Omega$},\\
\end{array}
\right.
\end{eqnarray*}
was considered, the author proved that there exist a positive, a
negative and a sigh-changing solution when the parameter $\lambda$
is greater than the second eigenvalue of the Steklov problem (\ref{e1.9}),
$f,g$ satisfying $\displaystyle\lim_{s\rightarrow
0}\frac{f(x,s)}{|s|^{p-2}s} = \displaystyle\lim_{s\rightarrow
0}\frac{g(x,s)}{|s|^{p-2}s} = 0$ and there exist $\delta_{f} > 0$
such that $\frac{f(x,s)}{|s|^{p-2}s} \geq 0$ when $0 < |s| <
\delta_{f}$, (It was proved in \cite{L} that the first eigenvalue of
the Steklov problem is isolated, so the second eigenvalue is the the
minimal eigenvalue greater than the first one). In \cite{AR09}, the
authors considered the following problem
\begin{eqnarray*}
\left \{
\begin{array}
{ll} -\Delta_{p}u + m(x)|u|^{p-2}u = \lambda a(x)|u|^{q-2}u,
&\textrm{  in $\Omega$},\\
|\nabla u|^{p-2}\frac{\partial u}{\partial n} = b(x)|u|^{r-2}u,
&\textrm{  on $\partial\Omega$},\\
\end{array}
\right.
\end{eqnarray*}
where  $1 < q < p < r < p^{*}$, $\|m\|_{\infty} > 0$, $a(x)\in
C(\overline{\Omega}),\|a\|_{\infty} = 1$ and $b(x)\in
C(\partial\Omega)$, $\|b\|_{\infty} = 1$. They proved that for $0 <
\lambda < \lambda^{*}$($\lambda^{*}$ is a constant depends on
$p,q,r$ and the best Sobolev constants of the embedding
$W^{1,p}_{0}(\Omega)\hookrightarrow L^{q}(\Omega)$ and
$W^{1,p}_{0}(\Omega)\hookrightarrow L^{r}(\partial\Omega)$), the
above problem has two solutions.

We note that all the problems listed above deal with the existence
or multiplicity problems for definite weight (i.e.the weight does
not change sign) or a restricted $\lambda$. However, Theorem 1.1-1.3, 1.1$'$-1.3$'$
and corollary 1.4-1.5 are for indefinite weight and every
$\lambda\in\mathbf{R}$.

For the no-flux problem (\ref{e1.6}), if we set $N=1$ and $\Omega=(0,T)$,
we get the following periodic problem for one-dimensional
$p$-Laplace equation:
\begin{eqnarray*}
\left \{
\begin{array}
{ll} -(|u'|^{p-2}u')' =\lambda V(x)|u|^{p-2}u + h(x,u),\\
u(0)=u(T),\\
u'(0)=u'(T).
\end{array}
\right.
\end{eqnarray*}
The periodic solution of $p$-laplace equation has been considered in
many papers, for example, \cite{BM07,BM08,MM00}. To the author's
knowledge, when applied to this one-dimensional case, our results as
stated in Theorem 1.3 and Corollary 1.4 are also new.

This paper is organized as follows. In section 2, we recall some
notations, definitions and some useful lemmas. In section 3, we
study the eigenvalue problems with Steklov, No-flux, Neumann, Robin
boundary value conditions respectively. We prove the existence of a
divergent sequence of eigenvalues  by critical point theory for even
functionals on Finsler manifolds. In section 4, we prove Theorems 1.1, 1.2
and Theorems 1.1$'$, 1.2$'$.

\section{Notations, definitions and known results}
Let $X$ be a closed linear subspace of $W^{1,p}(\Omega)$ such that
$W^{1,p}_{0}(\Omega) \subseteq X \subseteq W^{1,p}(\Omega)$ with the
norm $\|\cdot\|$ induced from the usual norm in $W^{1,p}(\Omega)$.
In this paper, we will also use an equivalent norm on $X$ defined by
$\|u\|_{\varepsilon}^{p} = \int_{\Omega}(|\nabla u|^{p} +
\varepsilon |u|^{p}){\rm d}x$ for a positive number $\varepsilon$.
By Pettis's theorem, $X$ is reflexive.
\subsection{Sobolev embedding theorem}
In the following, we will use Sobolev embedding theorem and trace
theorem frequently. So we list them as the following lemmas (see
\cite{KJF}).
\begin{lem}\label{l2.1} Let $\Omega$ be a bounded domain in
$\mathbf{R}^{N}$ with smooth boundary, there hold\\
(i)If $p < N$, then $W^{1,p}(\Omega) \hookrightarrow L^{q}(\Omega)$
for $1 \leq q \leq \frac{Np}{N-p}$, moreover, $W^{1,p}(\Omega)
\hookrightarrow\hookrightarrow
L^{q}(\Omega)$\\$~~~~$ when $1 \leq q < \frac{Np}{N-p} $,\\
(ii) If $p = N$, then $W^{1,p}(\Omega)
\hookrightarrow\hookrightarrow L^{q}(\Omega)$for $1 \leq q < \infty$,\\
(iii) If $p > N$, then $W^{1,p}(\Omega) \hookrightarrow
C^{1-\frac{N}{p}}(\overline{\Omega})$ and $W^{1,p}(\Omega)
\hookrightarrow\hookrightarrow C^{\beta}(\overline{\Omega})$ for $0
\leq \beta < 1-\frac{N}{p}$,\\
here and in the sequel, $\hookrightarrow$ means continuous embedding
map, and $\hookrightarrow\hookrightarrow$ means compact embedding
map.
\end{lem}
\begin{lem}\label{l2.2} Let $\Omega$ be a bounded domain in
$\mathbf{R}^{N}$ with smooth boundary, there hold\\
(i) If $p < N$, then $W^{1,p}(\Omega) \hookrightarrow
L^{q}(\partial\Omega)$ for $1 \leq q \leq \frac{Np-p}{N-p}$,
$W^{1,p}(\Omega) \hookrightarrow\hookrightarrow
L^{q}(\partial\Omega)$ for \\$~~~~$ $1 \leq q < \frac{Np-p}{N-p}$,\\
(ii) If $p = N$, then $W^{1,p}(\Omega)
\hookrightarrow\hookrightarrow L^{q}(\partial\Omega)$ for $1 \leq q < \infty$,\\
(iii) If $p > N$, we have $W^{1,p}(\Omega)
\hookrightarrow\hookrightarrow L^{q}(\partial\Omega)$ for $q \geq
1$.
\end{lem}
\subsection{Weak solution}
We give the following definitions on weak solution (See, for
example, \cite{L}, for details). \\(i) Let $u\in W^{1,p}(\Omega)$,
we say it is a weak solution of (\ref{e1.1}) if it satisfies the equation
\begin{eqnarray*}
\int_{\Omega}|\nabla u|^{p-2}\nabla u \cdot \nabla v {\rm d}x +
\int_{\Omega}\varepsilon|u|^{p-2}u v {\rm d}x = \lambda
\int_{\partial\Omega}V(x)|u|^{p-2}u v{\rm d}S_x +
\int_{\partial\Omega}h(x,u)v{\rm d}S_x
\end{eqnarray*}
~~~~~for any $v\in W^{1,p}(\Omega)$, \\
(ii) Let $u\in W^{1,p}_{0}(\Omega)\oplus \mathbf{R}$, we say it is a
weak solution of (\ref{e1.2}) if it satisfies the equation
\begin{eqnarray*}
\int_{\Omega}|\nabla u|^{p-2}\nabla u\cdot \nabla v{\rm d}x
+\int_{\Omega}\varepsilon|u|^{p-2}uv{\rm d}x=
\lambda\int_{\Omega}V(x)|u|^{p-2}uv{\rm d}x +
\int_{\Omega}h(x,u)v{\rm d}x
\end{eqnarray*}~~~~~for any $v\in W^{1,p}_{0}(\Omega)\oplus
\mathbf{R}$,\\
(iii) Let $u\in W^{1,p}(\Omega)$, we say it is a weak solution of
(\ref{e1.3}) if it satisfies the equation
\begin{eqnarray*}
\int_{\Omega}|\nabla u|^{p-2}\nabla u\cdot \nabla v{\rm d}x
+\int_{\Omega}\varepsilon|u|^{p-2}uv{\rm d}x=
\lambda\int_{\Omega}V(x)|u|^{p-2}uv{\rm d}x +
\int_{\Omega}h(x,u)v{\rm d}x
\end{eqnarray*}~~~~~~for any $v\in W^{1,p}(\Omega)$,\\
(iv) Let $u\in W^{1,p}(\Omega)$, we say it is a weak solution of the
(\ref{e1.4}) if it satisfies the equation
\begin{eqnarray*}
&&\int_{\Omega}|\nabla u|^{p-2}\nabla u\cdot \nabla v {\rm d}x +
\int_{\Omega}\varepsilon|u|^{p-2}uv{\rm d}x +
\int_{\partial\Omega}\gamma(x) |u|^{p-2}uv {\rm d}S_x = \lambda
\int_{\Omega}V(x)|u|^{p-2}uv {\rm d}x \\&&+ \int_{\Omega}h(x,u)v{\rm
d}x {\rm~for ~any~} v\in W^{1,p}(\Omega).\end{eqnarray*}
\subsection{Cohomological index}
In this subsection, we recall the construction and some properties
of the cohomological index of Fadell-Rabinowitz for a
$\mathbf{Z_{2}}$-set, see \cite{FR77,FR78,PAR} for details. For
simplicity, we only consider the usual $\mathbf{Z_{2}}$-action on a
linear space, i.e., $\mathbf{Z_{2}}=\{1,-1\}$ and the action is the
usual multiplication. In this case, the $\mathbf{Z_{2}}$-set $A$ is
a center symmetric set with $-A=A$.

Let $W$ be a normed linear space. We denote by $\mathcal {S}(W)$ the
set of all center symmetric subset of  $W$  not containing the
origin in $W$. For $A\in \mathcal {S}(W)$, denote $\bar{A} =
A/\mathbf{Z_{2}}$. Let $f: \bar{A} \rightarrow \mathbf{R}P^{\infty}$
be the classifying map and $f^{*}:
H^{*}(\mathbf{R}P^{\infty})=\mathbf{Z_{2}}[\omega] \rightarrow
H^{*}(\bar{A})$ the induced homomorphism of the cohomology rings.
The cohomological index of $A$, denoted by $i(A)$, is defined by
$\displaystyle\sup\{k \geq 1: f^{*}(\omega^{k-1}) \neq 0\}$. Here,
we list some properties which will be useful for us in this paper.
Let $A,B\in \mathcal {S}(W)$\\
(i1) (monotonicity) if $A \subseteq B$, then $i(A) \leq i(B)$.\\
(i2) (invariance) if $f:A \rightarrow B$ is an odd homeomorphism, then $i(A)=i(B)$. \\
(i3) (continuity) if $C$ is a closed symmetric subset of $A$, then
there exists a closed symmetric $~~~~~~$neighborhood $N$ of $C$ in
$A$, such that $i(N) = i(C)$, hence the interior of $N$ is also a
$~~~~~~$neighborhood of $C$ in $A$ and
$i({\rm int} N)=i(C)$.\\
(i4) (neighborhood of zero) if $U$ is bounded closed symmetric
neighborhood of the origin\\$~~~~~~$  in $W$, then $i(\partial U) =
\dim W$.\\ For more properties about the cohomological index, we
refer to \cite{PAR}.
\subsection{Some useful lemmas}
In this subsection, we recall some known results which will be
useful in section 3 and section 4. The first one is a linking
theorem for cones in normed spaces which is the theoretical tool of
this paper. It is contained in Corollary 2.9, Theorem 2.8,
Proposition 2.4 and Theorem 2.2 of \cite{DL07}. Here, we write it as
one lemma.
\begin{lem}\label{l2.3}(\cite{DL07})
Let $X$ be a real normed space and let $C_{-}$, $C_{+}$ be two
symmetric cones in $X$ such that $C_{+}$ is closed in $X$, $C_{-}
\cap C_{+} = \{0\}$ and  $$i(C_{-}\setminus \{0\}) = i(X\setminus
C_{+}) = m < \infty.$$ Define the following four sets by
\begin{eqnarray*}
&&D_{-}=\{u\in C_{-}:\|u\| \leq r_{-}\},\\
&&S_{+}=\{u\in C_{+}:\|u\|=r_{+}\},\\
&&Q=\{u + te: u\in C_{-}, t \geq 0,\|u +te\| \leq r_{-}\},\;\; e\in X\setminus C_-,\\
&&H=\{u + te:u\in C_{-}, t \geq 0, \|u + te\| = r_{-}\}.
\end{eqnarray*}
Then $(Q, D_{-} \cup H)$ links $S_{+}$ cohomologically in dimension
$m+1$ over $\mathbf{Z_{2}}$. Moreover, suppose $f\in C^1(X,
\mathbf{R})$ satisfying the $(PS)$ condition, and
$\displaystyle\displaystyle\sup_{x\in D_{-} \cup H}f(x) <
\displaystyle\inf_{x\in S^{+}}f(x)$,
$\displaystyle\displaystyle\sup_{x\in Q}f(x) < \infty$. Then $f$ has
a critical value $c \geq \displaystyle\displaystyle\inf_{x\in
S^{+}}f(x)$.
\end{lem}

{\bf Remark:} Recently, in \cite{D09},
the author extended it to more general case
(the functional space is completely regular topological space or metric space).
If the functional space $X$ is a real Banach space,
according to the proof of Theorem 6.10 in \cite{D09},
the Cerami condition is sufficient for the compactness of the set of critical points at a fixed level and the first deformation lemma to hold (see \cite{PAR}).
So this critical point theorem still hold under the Cerami condition.

The results in section 3 is based on the following theorem.
\begin{lem}\label{l2.4}(Proposition $3.52$ in \cite{PAR})
Suppose $\mathcal{M}$ is a $C^{1}$ Finsler manifold with free
$\mathbf{Z}_{2}$-action, $\Phi\in C^{1}(\mathcal{M}, \mathbf{R})$
and $\Phi$ is even (i.e.$\mathbf{Z}_{2}$-invariant). Set
$$\mathcal{F} _{k}=\{M:M ~{\rm is}~\mathbf{Z}_{2}{\rm-invariant~and} ~i(M)
\geq k\}\;\;{\rm and}\;\; c_{k}=\displaystyle\inf_{M\in
\mathcal{F}_{k}}\displaystyle\sup_{u\in M}\Phi(u).$$
Then the following two statements are true:\\
(i) If $-\infty < c_{k} = \cdots=c_{k+m-1}=c < +\infty$ and $\Phi$
satisfies $(PS)_{c}$, then we have $i(K^{c}) \geq m$.
\\$~~~~$Moreover, if $-\infty < c_{k} \le \cdots\le c_{k+m-1} <
+\infty$ and the functional $\Phi$ satisfies $(PS)_{c}$ for $~~~$ $c
= c_{k},\cdots,c_{k+m-1}$, then all $c_{k},\cdots,c_{k+m-1}$ are
critical values and $\Phi$  has at least $m$\\$~~~$ distinct pairs
of
critical points.\\
(ii) If $-\infty < c_{k} < +\infty$ for all sufficiently large $k$
and $\Phi$ satisfies $(PS)$, then $c_{k}\nearrow +\infty$.
\end{lem}
In the proof of the main results, we will also use the following
technical lemma.
\begin{lem}\label{l2.5}(Lemma $4.2$ in \cite{DL07})
Let $E$ be a measurable subset of $\mathbf{R}^{n}$, let $1 \leq
\alpha < \infty$, $1 \leq \beta < \infty$ and $h : E \times
\mathbf{R} \rightarrow \mathbf{R}$ be Carath\'{e}odory function.
Assume that, for every $\epsilon > 0$, there exists $a_{\epsilon}\in
L^{\beta}(E)$ such that $|h(x,s)| \leq a_{\epsilon}(x) + \epsilon
|s|^{\frac{\alpha}{\beta}}$ for a.e.$x\in E$ and every
$s\in\mathbf{R}$. Then, if $(u_{k})$ is a bounded sequence  in
$L^{\alpha}(E)$ and convergent to $u$ a.e.in $E$, we have that
$(h(x,u_{k}))$ is convergent to $h(x,u)$ strongly in $L^{\beta}(E)$.
\end{lem}
\begin{rem}\label{l2.6}If the condition $|h(x,s)| \leq a_{\epsilon}(x) + \epsilon
|s|^{\frac{\alpha}{\beta}}$ only holds for a.e.$x\in E$ and $|s|
\leq S$($S$ is a positive constant), the conclusion also holds if
$\|u_{n}\|_{\infty} \leq S$, $\|u\|_{\infty} \leq S$.\end{rem}
\section{Existence of a divergent sequence of eigenvalues}
In this section, we assume that $meas\{x\in \Omega: V(x) > 0\} > 0$
if $V$ is defined on $\Omega$, $meas\{x\in \partial\Omega: V(x) >
0\}> 0$ if $V$ is defined on $\partial\Omega$. We consider the following eigenvalue problems\\
Steklov problem
\begin{eqnarray*}S(\Omega)_{\varepsilon}\quad\;\;
\left \{
\begin{array}
{ll}\Delta_{p}u = \varepsilon|u|^{p-2}u,
&\textrm{  in $\Omega$},\\
|\nabla u|^{p-2}\frac{\partial u}{\partial n}=\lambda V(x
)|u|^{p-2}u,
&\textrm{  on $\partial\Omega$},\\
\end{array}
\right.
\end{eqnarray*}
No-flux problem
\begin{eqnarray*}P(\Omega)_{\varepsilon}\quad\;\;
\left \{
\begin{array}
{ll} -\Delta_{p}u + \varepsilon |u|^{p-2}u=\lambda V(x)|u|^{p-2}u ,
&\textrm{  in $\Omega$},\\
u=constant,
&\textrm{  on $\partial\Omega$},\\
\int_{\partial\Omega}|\nabla u|^{p-2}\frac{\partial u}{\partial
n}{\rm d}S_x=0,
\end{array}
\right.
\end{eqnarray*}
Neumann problem
\begin{eqnarray*}N(\Omega)_{\varepsilon}\quad\;\;
\left \{
\begin{array}
{ll} -\Delta_{p}u +\varepsilon |u|^{p-2}u=\lambda V(x)|u|^{p-2}u,
&\textrm{  in $\Omega$},\\
\frac{\partial u}{\partial n}=0,
&\textrm{  on $\partial\Omega$},\\
\end{array}
\right.
\end{eqnarray*}
Robin problem
\begin{eqnarray*}R(\Omega)_{\varepsilon}\quad\;\;
\left \{
\begin{array}
{ll} -\Delta_{p}u +\varepsilon |u|^{p-2}u=\lambda V(x)|u|^{p-2}u,
&\textrm{  in $\Omega$},\\
|\nabla u|^{p-2}\frac{\partial u}{\partial n}+\gamma(x)
|u|^{p-2}u=0,
&\textrm{  on $\partial\Omega$}.\\
\end{array}
\right.
\end{eqnarray*}
In \cite{L}, for $V \equiv 1$ and $\varepsilon =1$ in the case
$S(\Omega)_{\varepsilon}$, $\varepsilon = 0$ in the cases
$P(\Omega)_{\varepsilon}$, $N(\Omega)_{\varepsilon}$,
$R(\Omega)_{\varepsilon}$, the author proved that these four
problems has a divergent sequence of eigenvalues respectively by
Ljusternik-Schnirelman principle. By critical point theory for
functionals on Finsler manifolds, we also get a divergent sequence
of eigenvalues respectively.
\subsection{A general eigenvalue problem}
Let $a\in L^{r}(\Omega)$, $b\in L^{r}(\partial\Omega)$, $\beta \in
L^{\infty}(\partial\Omega)$ and $\beta(x) \geq 0$ for a.e.$x\in
\partial\Omega$. We suppose that $a$, $b$ satisfy the following
assumption:\\
({\bf A}): If $meas\{x\in \Omega: a(x) > 0\} = 0$, then $a \equiv
0$, $meas\{x\in \partial\Omega: b(x) > 0\} > 0$ and
$~~~~~~~~X=W^{1,p}(\Omega)$.

Define on $X$ the functional
\begin{eqnarray*}
F(u) = \frac{1}{p}\int_{\Omega}a(x)|u(x)|^{p}{\rm d}x +
\frac{1}{p}\int_{\partial\Omega}b(s)|u(x)|^{p}{\rm d}S_x,
\end{eqnarray*}
\begin{eqnarray*}
G_{\varepsilon}(u) = \frac{1}{p}\int_{\Omega}(|\nabla
u|^{p}+\varepsilon|u|^{p}){\rm d}x+\frac{1}{p}
\int_{\partial\Omega}\beta(s)|u(x)|^{p}{\rm d}S_x.
\end{eqnarray*}
We want to solve the problem
\begin{equation}\label{e3.1}
G'_{\varepsilon}(u)=\lambda F'(u).
\end{equation}
Clearly, we have \begin{eqnarray*}F\in C^{1}, ~\langle
F'(u),v\rangle = \int_{\Omega}a|u|^{p-2}uv{\rm d}x +
\int_{\partial\Omega}b|u|^{p-2}uv{\rm d}S_x,\end{eqnarray*} and
\begin{eqnarray*}G_{\varepsilon}\in C^{1}, ~\langle G'_{\varepsilon}(u),v\rangle
= \int_{\Omega}(\nabla u|^{p-2}\nabla u \cdot \nabla v + \varepsilon
|u|^{p-2}uv){\rm d}x + \int_{\partial\Omega}\beta |u|^{p-2}uv{\rm
d}S_x.\end{eqnarray*}

First, we consider the case $\varepsilon > 0$.
\begin{lem}\label{l3.1} For any $u,v\in X$,
we have \begin{eqnarray*}\langle
G'_{\varepsilon}(u)-G'_{\varepsilon}(v),u-v\rangle \geq
(\|u\|^{p-1}_{\varepsilon}-\|v\|^{p-1}_{\varepsilon})(\|u\|_{\varepsilon}-\|v\|_{\varepsilon})\end{eqnarray*}
\end{lem}
\pf: Its proof is the same as Lemma 2.3 in \cite{L}. For reader's
convenience we give it here.

By direct computations, we have
\begin{eqnarray*}
\langle G'_{\varepsilon}(u)-G'_{\varepsilon}(v),u-v\rangle &=&
\int_{\Omega}[|\nabla u|^{p}+|\nabla v|^{p}-|\nabla u|^{p-2}\nabla
u\cdot \nabla v-|\nabla
v|^{p-2}\nabla v\cdot \nabla u]dx\\
&+&\varepsilon\int_{\Omega}
(|u|^{p}+|v|^{p}-|u|^{p-2}uv-|v|^{p-2}vu)dx\\
&+&\int_{\partial\Omega}
\beta(|u|^{p}+|v|^{p}-|u|^{p-2}uv-|v|^{p-2}vu){\rm d}S_{x}.
\end{eqnarray*}
It follows from the proof of Lemma $2.3$ in \cite{L} that
\begin{eqnarray*}
\int_{\partial\Omega}
\beta(|u|^{p}+|v|^{p}-|u|^{p-2}uv-|v|^{p-2}vu){\rm d}S_{x} \geq 0.
\end{eqnarray*}
Hence
\begin{eqnarray*} \langle
G'_{\varepsilon}(u)-G'_{\varepsilon}(v),u-v\rangle &\geq &
\int_{\Omega}[|\nabla u|^{p}+|\nabla v|^{p}-|\nabla u|^{p-2}\nabla
u\cdot \nabla v-|\nabla v|^{p-2}\nabla v\cdot \nabla
u]dx \\
&+&
\varepsilon\int_{\Omega}(|u|^{p}+|v|^{p}-|u|^{p-2}uv-|v|^{p-2}vu)dx\\
&=& \|u\|_{\varepsilon}^{p}+\|v\|_{\varepsilon}^{p}
-\int_{\Omega}(|\nabla u|^{p-2}\nabla u\cdot \nabla v+\varepsilon
|u|^{p-2}uv)dx \\
&-&\int_{\Omega}(|\nabla v|^{p-2}\nabla v\cdot \nabla u+\varepsilon
|v|^{p-2}vu)dx.
\end{eqnarray*}

Applying H\"{o}lder inequality, we have
\begin{eqnarray*}
&&\int_{\Omega}(|\nabla u|^{p-2}\nabla u\cdot \nabla v+\varepsilon
|u|^{p-2}uv)dx
\\&&\leq\Big(\int_{\Omega}|\nabla
u|^{p}dx\Big)^{\frac{p-1}{p}}\Big(\int_{\Omega}|\nabla
v|^{p}dx\Big)^{\frac{1}{p}}+\Big(\int_{\Omega}\varepsilon
|u|^{p}dx\Big)^{\frac{p-1}{p}}\Big(\int_{\Omega}\varepsilon
|v|^{p}dx\Big)^{\frac{1}{p}}.
\end{eqnarray*}
Similar to the proof of Lemma 2.3 in \cite{L}, we use the following
inequality
\begin{eqnarray*}(a+b)^{\alpha}(c+d)^{1-\alpha} \geq
a^{\alpha}c^{1-\alpha}+b^{\alpha}d^{1-\alpha}
\end{eqnarray*}
which holds for any $\alpha \in (0,1)$ and for any $a > 0$, $b > 0$,
$c
> 0$, $d
> 0 $. Set \begin{eqnarray*}
a=\int_{\Omega}|\nabla u|^{p}dx,  b=\int_{\Omega}\varepsilon
|u|^{p}dx,  c=\int_{\Omega}|\nabla v|^{p}dx,
d=\int_{\Omega}\varepsilon |v|^{p}dx, \alpha = \frac{p-1}{p},
\end{eqnarray*}
we can deduce that
\begin{eqnarray*} \int_{\Omega}(|\nabla
u|^{p-2}\nabla u\cdot \nabla v+\varepsilon |u|^{p-2}uv)dx \leq
\|u\|_{\varepsilon}^{p-1}\|v\|_{\varepsilon}.
\end{eqnarray*}

Similarly, we can obtain
\begin{eqnarray*}
\int_{\Omega}(|\nabla v|^{p-2}\nabla v\cdot \nabla u+\varepsilon
|v|^{p-2}vu)dx \leq \|v\|_{\varepsilon}^{p-1}\|u\|_{\varepsilon}.
\end{eqnarray*}

Therefore, we have  \begin{eqnarray*} \langle
G'_{\varepsilon}(u)-G'_{\varepsilon}(v),u-v\rangle &\geq &
\|u\|_{\varepsilon}^{p} + \|v\|_{\varepsilon}^{p} -
\|u\|_{\varepsilon}^{p-1}\|v\|_{\varepsilon} -
\|v\|_{\varepsilon}^{p-1}\|u\|_{\varepsilon}\\
&= & (\|u\|_{\varepsilon}^{p-1} - \|v\|_{\varepsilon}^{p-1})
(\|u\|_{\varepsilon} - \|v\|_{\varepsilon})\\
&\geq & 0.
\end{eqnarray*} \qed
\begin{lem}\label{l3.2} If $u_{n}\rightharpoonup u$, $\langle
G'_{\varepsilon}(u_{n}),u_{n}-u\rangle \rightarrow 0$, then
$u_{n}\rightarrow u$ in $X$.
\end{lem}
\pf: By Sobolev's compact embedding theorem we have $u_{n}
\rightarrow u$ in $L^{p}(\Omega)$. Since $X$ is a reflexive Banach
space, weak convergence and norm convergence imply strong
convergence(see the proof of Proposition 2.4 in \cite{L}). So we
only need to show that $\|u_{n}\|_{\varepsilon} \rightarrow
\|u\|_{\varepsilon}$.

Notice that
\begin{eqnarray*}
\lim_{n\rightarrow \infty}\langle
G'_{\varepsilon}(u_{n})-G'_{\varepsilon}(u),u_{n}-u\rangle
=\lim_{n\rightarrow \infty}(\langle
G'_{\varepsilon}(u_{n}),u_{n}-u\rangle - \langle
G'_{\varepsilon}(u),u_{n}-u\rangle) = 0.
\end{eqnarray*}
By the Lemma \ref{l3.1} we have
\begin{eqnarray*}
\langle G'_{\varepsilon}(u_{n})-G'_{\varepsilon}(u), u_{n}-u\rangle
\geq
(\|u\|_{\varepsilon}^{p-1}-\|u\|_{\varepsilon}^{p-1})(\|u_{n}\|_{\varepsilon}-\|u\|_{\varepsilon})
\geq 0.
\end{eqnarray*}
Hence $\|u_{n}\|_{\varepsilon} \rightarrow \|u\|_{\varepsilon}$ as
$n \rightarrow \infty$ and the assertion follows.\qed
\begin{lem}\label{l3.3} $F^{'}$ is weak-to-strong continuous, i.e.
$u_{n}\rightharpoonup u$ in $X$ implies $F'(u_{n})\rightarrow
F'(u)$.
\end{lem}
\pf: Let $u_{n}\rightharpoonup u$ in $X$. We have to show that
$F'(u_{n})\rightarrow F'(u)$ in $X^{*}$. The proof is similar to the
proof of Proposition 2.2 in \cite{L}.
\par
{\bf If $1 < p < N$.} For any $v\in X$, by H\"{o}lder inequality ,
Sobolev embedding theorem and the identity
\begin{eqnarray*}
\frac{p}{N} + \frac{p-1}{\frac{Np}{N-p}} + \frac{1}{\frac{Np}{N-p}}
= 1, \frac{p-1}{N-1} + \frac{p-1}{\frac{Np-p}{N-p}} +
\frac{1}{\frac{Np-p}{N-p}} = 1,
\end{eqnarray*}
we have that
\begin{eqnarray*}
&&|\langle F'(u_{n})-F'(u),v\rangle|\\&&\leq
|\int_{\Omega}a(|u_{n}|^{p-2}u_{n}-|u|^{p-2}u)v{\rm d}x|+
|\int_{\partial\Omega}b(|u_{n}|^{p-2}u_{n}-|u|^{p-2}u)v{\rm d}S_x|\\
&&\leq
C_{1}^{'}\|a\|_{L^{r}(\Omega)}\||u_{n}|^{p-2}u_{n}-|u|^{p-2}u\|_{L^{\frac{\beta}{p-1}}(\Omega)}\|v\|_{L^{\frac{Np}{N-p}}(\Omega)}\\
&&+
C_{2}^{'}\|b\|_{L^{r}(\partial\Omega)}\||u_{n}|^{p-2}u_{n}-|u|^{p-2}u\|_{L^{\frac{\gamma}{p-1}}(\partial\Omega)}
\|v\|_{L^{\frac{Np-p}{N-p}}(\partial\Omega)}\\
&&\leq
C_{1}\|a\|_{L^{r}(\Omega)}\||u_{n}|^{p-2}u_{n}-|u|^{p-2}u\|_{L^{\frac{\beta}{p-1}}(\Omega)}\|v\|\\
&&+
C_{2}\|b\|_{L^{r}(\partial\Omega)}\||u_{n}|^{p-2}u_{n}-|u|^{p-2}u\|_{L^{\frac{\gamma}{p-1}}
(\partial\Omega)}\|v\|.
\end{eqnarray*} Here $\beta$ and $\gamma$ satisfy $\max\{p - 1,1\} < \beta
< \frac{Np}{N-p}$, and $\max\{p -1,1\} < \gamma < \frac{Np-p}{N-p}$.

To prove the conclusion, we only need to show that
$|u_{n}|^{p-2}u_{n} \rightarrow |u|^{p-2}u$ in
$L^{\frac{\beta}{p-1}}(\Omega)$ and $|u_{n}|^{p-2}u_{n} \rightarrow
|u|^{p-2}u$ in $L^{\frac{\gamma}{p-1}}(\partial\Omega)$. To see
this, let $w_{n}=|u_{n}|^{p-2}u_{n}$ and $w=|u|^{p-2}u$. Since
$u_{n}\rightharpoonup u$ in $W^{1,p}(\Omega)$, $u_{n} \rightarrow u$
in $L^{\beta}(\Omega)$, it follows that  $w_{n}(x) \rightarrow
w(x)$, a.e. in $\Omega$ and
$\int_{\Omega}|w_{n}|^{\frac{\beta}{p-1}}{\rm d}x \rightarrow
\int_{\Omega}|w|^{\frac{\beta}{p-1}}{\rm d}x$, by Proposition $2.4$
in \cite{FZ03}, we conclude that $w_{n} \rightarrow w$ in
$L^{\frac{\beta}{p-1}}(\Omega)$. The proof of $u_{n} \rightarrow u$
in $L^{\frac{\gamma}{p-1}}(\partial\Omega)$ is similar.

{\bf If $ p > N$.} For any $v\in X$, we have
\begin{eqnarray*}
&&|\langle F'(u_{n})-F'(u),v\rangle|\\&&\leq
|\int_{\Omega}a(|u_{n}|^{p-2}u_{n}-|u|^{p-2}u)v{\rm d}S_x| +
|\int_{\partial\Omega}b(|u_{n}|^{p-2}u_{n}-|u|^{p-2}u)v{\rm d}S_x|\\
&&\leq
\|a\|_{L^{1}(\Omega)}\||u_{n}|^{p-2}u_{n}-|u|^{p-2}u\|_{L^{\infty}(\Omega)}\|v\|_{L^{\infty}(\Omega)}\\
&&+\|b\|_{L^{1}(\partial\Omega)}\||u_{n}|^{p-2}u_{n}-|u|^{p-2}u\|_{L^{\infty}(\partial\Omega)}
\|v\|_{L^{\infty}(\partial\Omega)}\\
&&\leq
C_{1}\|a\|_{L^{1}(\Omega)}\||u_{n}|^{p-2}u_{n}-|u|^{p-2}u\|_{L^{\infty}(\Omega)}\|v\|\\
&&+
C_{2}\|b\|_{L^{1}(\partial\Omega)}\||u_{n}|^{p-2}u_{n}-|u|^{p-2}u\|_{L^{\infty}
(\partial\Omega)}\|v\|.
\end{eqnarray*}

By the Sobolev embedding theorem, we have that $u_{n},u\in
C(\bar{\Omega})$ and $u_{n} \rightarrow u$ uniformly, so
$\||u_{n}|^{p-2}u_{n} - |u|^{p-2}u\|_{L^{\infty}(\Omega)}
\rightarrow 0$ and $\||u_{n}|^{p-2}u_{n} -
|u|^{p-2}u\|_{L^{\infty}(\partial\Omega)} \rightarrow 0$. So the
conclusion follows in this case.

{\bf If $p = N$.}  For any $v\in X$, by H\"{o}lder inequality and
the Sobolev embedding theorem it follows that
\begin{eqnarray*}
&&|\langle F'(u_{n})-F'(u),v\rangle|\\&&\leq
|\int_{\Omega}a(|u_{n}|^{p-2}u_{n}-|u|^{p-2}u)v{\rm d}x| +
|\int_{\partial\Omega}b(|u_{n}|^{p-2}u_{n}-|u|^{p-2}u)v{\rm d}S_x|\\
&&\leq
\|a\|_{L^{r}(\Omega)}\||u_{n}|^{p-2}u_{n}-|u|^{p-2}u\|_{L^{\frac{\beta}{p-1}}(\Omega)}\|v\|_{L^{s}(\Omega)}\\
&&+\|b\|_{L^{r}(\partial\Omega)}\||u_{n}|^{p-2}u_{n}-|u|^{p-2}u\|_{L^{\frac{\gamma}{p-1}}(\partial\Omega)}
\|v\|_{L^{t}(\partial\Omega)}\\ &&\leq
C_{1}\|a\|_{L^{r}(\Omega)}\||u_{n}|^{p-2}u_{n}-|u|^{p-2}u\|_{L^{\frac{\beta}{p-1}}(\Omega)}\|v\|\\
&&+C_{2}\|b\|_{L^{r}(\partial\Omega)}\||u_{n}|^{p-2}u_{n}-|u|^{p-2}u\|_{L^{\frac{\gamma}{p-1}}
(\partial\Omega)}\|v\|.
\end{eqnarray*}
Here $\beta$ and $\gamma$ satisfy $\beta > \max\{p - 1,1\}$, $\gamma
> \max\{p -1,1\}$, and $s,t
> 1$ are real number such that \begin{eqnarray*}\frac{1}{r} + \frac{p-1}{\beta} + \frac{1}{s} =
1,\frac{1}{r} + \frac{p-1}{\gamma} + \frac{1}{t} = 1.\end{eqnarray*}

To prove the conclusion, we only need to show that
$|u_{n}|^{p-2}u_{n} \rightarrow |u|^{p-2}u$ in
$L^{\frac{\beta}{p-1}}(\Omega)$ and $|u_{n}|^{p-2}u_{n} \rightarrow
|u|^{p-2}u$ in $L^{\frac{\gamma}{p-1}}(\partial\Omega)$. The proof
is similar to the case $p < N$.\qed
\begin{lem}\label{l3.4}
If $u_{n}\rightharpoonup u$, then $F(u_{n}) \rightarrow F(u)$.
\end{lem}
\pf: By the definition of $F$, there holds
\begin{eqnarray*}
p|F(u_{n}) - F(u)| &=& |\langle
F'(u_{n}),u_{n}\rangle - \langle F'(u),u\rangle|\\ &=& |\langle
F'(u_{n}) - F'(u), u_{n}\rangle + \langle F'(u),
u_{n} - u\rangle|\\
&\leq & \|F'(u_{n}) - F'(u)\| \|u_{n}\| + o(1).
\end{eqnarray*}
Because $u_{n} \rightharpoonup u$, $u_{n}$ is bounded. From Lemma \ref{l3.3}
, we have $F(u_{n}) \rightarrow F(u)$.\qed

Set $\mathcal{M} = \{u\in X: F(u) = 1\}$, it is nonempty since
$meas\{x\in\Omega: a(x) > 0\} > 0$ or $meas\{x\in\Omega: a(x) > 0\}
= 0$, $meas\{x\in \partial\Omega: b(x) > 0\} > 0$ and
$X=W^{1,p}(\Omega)$ by assumption({\bf A}) (for detail, see the
proof of Lemma \ref{l3.7}). Clearly $F(u)=\frac{1}{p}\langle
F^{'}(u),u\rangle$, so $1$ is a regular value of the function $F$.
Hence $\mathcal{M}$ is a $C^{1}$-Finsler manifold by the implicit
theorem. It is complete, symmetric, since $F$ is continuous and
even. Moreover, $0$ is not contained in $\mathcal{M}$, so the usual
$\mathbf{Z}_{2}$-action on $\mathcal{M}$ is free. Set
$\widetilde{G}_{\varepsilon} = G_{\varepsilon}|_{\mathcal{M}}$.
\begin{lem}\label{l3.5}
If $u\in\mathcal{M}$ satisfies $\widetilde{G}_{\varepsilon}(u) =
\lambda$ and $\widetilde{G}_{\varepsilon}'(u) = 0$, then
$(\lambda,u)$ is a solution to (\ref{e3.1})
\end{lem}
\pf: By Proposition $3.54$ in \cite{PAR}, the norm of
$\widetilde{G}_{\varepsilon}'(u)\in T^{*}_{u}\mathcal{M}$ is given
by $\|\widetilde{G}_{\varepsilon}'(u)\|_{u}^{*} =
\displaystyle\displaystyle\min_{\mu\in \mathbf{R}}\|G'(u) - \mu
F'(u)\|^{*}$(here the norm $\|\cdot\|^{*}_{u}$ is the norm in the
fibre $T^{*}_{u}\mathcal{M}$, and $\|\cdot\|^{*}$ is the operator
norm). Hence there exist $\mu\in \mathbf{R}$ such that
$G_{\varepsilon}'(u) - \mu F'(u) = 0$, that is $(\mu, u)$ is a
solution of (\ref{e3.1}) and $ \lambda=\widetilde{G}_{\varepsilon}(u) =
\mu$. \qed
\begin{lem}\label{l3.6} $\widetilde{G}_{\varepsilon}$ satisfies the
$(PS)$ condition, i.e. if $(u_{n})$ is a sequence on $\mathcal{M}$
such that $\widetilde{G}_{\varepsilon}(u_{n}) \rightarrow c$, and
$\widetilde{G}_{\varepsilon}'(u_{n}) \rightarrow 0$, then up to a
subsequence $u_{n} \rightarrow u\in \mathcal{M}$ in $X$
\end{lem}
\pf: First, from the definition of $G_{\varepsilon}$, we can deduce
that $(u_{n})$ is bounded. Since $X$ is reflexive, up to a
subsequence, $u_{n}$ converges weakly to some $u\in X$.

From $\widetilde{G}_{\varepsilon}'(u_{n}) \rightarrow 0$, we have
$G_{\varepsilon}'(u_{n}) - \mu_{n} F'(u_{n}) \rightarrow 0$ for a
sequence of real numbers $(\mu_{n})$. Then applying this formula to
$u_{n}$, we get $\mu_{n} \rightarrow c$. By Lemma \ref{l3.3},
$G_{\varepsilon}'(u_{n})\rightarrow c F'(u)$. Hence $\langle
G_{\varepsilon}'(u_{n}), u_{n} - u \rangle \rightarrow 0$. By Lemma
\ref{l3.2}, we get $u_{n} \rightarrow u$. \qed

Let $\mathcal{F}$ denote the class of symmetric subsets of
$\mathcal{M}$, let $\mathcal{F}_{n} = \{M\in \mathcal{F}: i(M) \geq
n\}$ and $\lambda_{n,\varepsilon} = \displaystyle\inf_{M\in
\mathcal{F}_{n}}\displaystyle\sup_{u\in
M}\widetilde{G}_{\varepsilon}(u)$. Since $\mathcal{F}_{n}
\displaystyle\supset \mathcal{F}_{n+1}$, $\lambda_{n,\varepsilon}
\leq \lambda_{n+1,\varepsilon}$.
\begin{lem}\label{l3.7}There exists a compact set in
$\mathcal{F}_{n}$.
\end{lem}
\pf: If $meas\{x\in \Omega: a(x) > 0\} = 0$, then by assuption({\bf
A}), $a\equiv 0$, $meas\{x\in \partial\Omega: b(x) > 0\} > 0$ and
$X=W^{1,p}(\Omega)$. We follow the idea in the proof of Theorem 3.2
in \cite{HT}. In this case, we can infer that $\forall n\in
\mathbf{N}^{*}$, there exist $n$ open balls $(B_{i})_{1 \leq i \leq
n}$ in $\partial\Omega$ such that $B_{i}\cap B_{j} = \emptyset$ if
$i \neq j$ and $meas(\{x\in
\partial\Omega: b(x) > 0\}\cap B_{i}) > 0$. Approximating the characteristic
function $\chi_{\{x\in \partial\Omega: b(x) > 0\}\cap B_{i}}$ by
$C^{\infty}(\partial\Omega)$ functions in
$L^{\frac{rp}{r-1}}(\partial\Omega)$, we can infer that there exists
a sequence $(u_{i})_{1 \leq i \leq n}\subseteq
C^{\infty}(\partial\Omega)$ such that
$\int_{\partial\Omega}b(s)|u_{i}|^{p}ds > 0 $ for all $i = 1,...n$
and supp$u_{i}\cap$ supp$u_{j}=\emptyset$ when $i \neq j$. From
trace theorem, we can find a sequence $(w_{i})_{1 \leq i \leq n}\in
X$ such that $\Gamma(w_{i}) = u_{i}$, here $\Gamma$ is the trace
map. So $F(w_{i})=\frac{1}{p}\int_{\partial\Omega}b(s)|u_{i}|^{p}ds>
0$. Normalizing $w_{i}$, we assume that $F(w_{i}) = 1$. Denote
$W_{n}$ the space generated by $(w_{i})_{1 \leq i \leq n}$. $\forall
w\in W_{n}$, we have $w = \sum_{i=1}^{n}\alpha_{i}w_{i}$ and
$F(w)=\sum_{i=1}^{n}|\alpha_{i}|^{p}$. So $w\to
\Big(F(w)\Big)^{\frac{1}{p}}$ defines a norm on $W_{n}$. Since
$W_{n}$ is finite-dimensional, this norm is equivalent to
$\|\cdot\|_{\varepsilon}$. So $\{w\in W_{n}: F(w)=1\}\subseteq
\mathcal{M}$ is compact with respect to the norm
$\|\cdot\|_{\varepsilon}$ and by $(i4)$ in section $2.3$, $i(\{w\in
W_{n}:F(w)=1\})=n$. So $\{w\in W_{n}: F(w)=1\}\in \mathcal{F}_{n}$.

If $meas\{x\in \Omega: a(x) > 0\} > 0$, the proof is similar, see
also the proof of Theorem 3.2 in \cite{HT}.\qed

Hence, $\lambda_{k,\varepsilon}$ is finite. Finally, from Lemma \ref{l2.4}
and Lemma \ref{l3.6}, we have $\lambda_{n,\varepsilon}$ is a divergent
sequence of critical values of $\widetilde{G}_{\varepsilon}$. So by
Lemma \ref{l3.5} we get a divergent sequence of eigenvalues for problem
(\ref{e3.1}).
\begin{lem}\label{l3.8} There holds
$$\lambda_{n,\varepsilon} = \displaystyle\inf_{K\in \mathcal{F}^{c}_{n}}\displaystyle\sup_{u\in
K}G_{\varepsilon}(u),$$ where $\mathcal{F}^{c}_{n}=\{K\in \mathcal
{F}_n: K ~is~ compact \}$.
\end{lem}
\pf: Indeed, the same reason as the proof of Proposition $3.1$ in
\cite{DL07}, we have that for every symmetric, open subset $A$ of
$\mathcal{M}$, $i(A) = \displaystyle\sup\{i(K): K ~is ~compact~ and
~ symmetric~
\\with ~ K \subseteq A\}$. This combines $(i3)$ in section $2.3$,
can deduce the assertion easily.\qed

Next, we consider the case $\varepsilon = 0$.

Put $G(u) = G_{0}(u) = \frac{1}{p}\int_{\Omega}|\nabla u|^{p}{\rm
d}x+\frac{1}{p}\int_{\partial\Omega}\beta(s)|u(s)|^{p}{\rm d}S_x$
and
$\lambda_{n}=\displaystyle\inf_{K\in \mathcal{F}^{c}_{n}}\displaystyle\sup_{u\in K}G(u)$. \\
To solve the eigenvalue problem $G^{'}(u) = \lambda F^{'}(u)$, we
follow the method in \cite{HT}.
\begin{lem}\label{l3.9} We have the following two statements:
\\
$(i)$ $\displaystyle\lim_{\varepsilon \rightarrow
0^{+}}\lambda_{n,\varepsilon}=\lambda_{n}$,\\
$(ii)$ $\lambda_{n} \rightarrow +\infty$ as $n \rightarrow +\infty$.
\end{lem}
\pf: (i) Let $\varepsilon > 0$, from the definition, we have
$\lambda_{n,\varepsilon} \geq \lambda_{n}$. $\forall \delta > 0$,
there exist $K = K(\delta)\in \mathcal{F}^{c}_{n}$ such that
$\lambda_{n} \leq \displaystyle\sup_{u\in K}G(u) < \lambda_{n} +
\delta$. Set $\gamma=\displaystyle\sup_{u\in K}\|u\|_{p}^{p}$, then there holds
$$\lambda_{n} \leq \lambda_{n,\varepsilon} \leq
\displaystyle\sup_{u\in K}G(u) + \frac{\varepsilon\gamma}{p}.$$ When
$\varepsilon$ is sufficiently small, we obtain
$\displaystyle\sup_{u\in K}G(u)+\frac{\varepsilon\gamma}{p}\leq
\lambda_{n} + \delta$. Thus $\lambda_{n} \leq
\lambda_{n,\varepsilon} \leq \lambda_{n} + \delta$ for all
$\varepsilon$ small enough. From this we get the desired result.

(ii) Fix $a(x)\in L^{r}(\Omega)$, $b(x)\in L^{r}(\partial\Omega)$,
since $\lambda_{n}$, $F$ and $\mathcal{F}_{c}^{n}$ depends on $a$
and $b$, we write $\lambda_{n}=\lambda_{n}(a,b)$ and $
F(u)=F(a,b)(u)$, $\mathcal{F}_{c}^{n}=\mathcal{F}_{c}^{n}(a,b)$. Let
$\tau
> 0$ be small, define
\begin{eqnarray*}\bar{a}(x) =
\left \{
\begin{array}
{ll} a(x),
&\textrm{  if $a(x) \geq \tau$},\\
\tau,
&\textrm{  if $a(x) < \tau$},\\
\end{array}
\right.
\end{eqnarray*}
and
\begin{eqnarray*}\bar{b}(x) =
\left \{
\begin{array}
{ll} b(x),
&\textrm{  if $b(x) \geq \tau$},\\
\tau,
&\textrm{  if $b(x) < \tau$}.\\
\end{array}
\right.
\end{eqnarray*}
Then $\bar{a}$, $\bar{b}$ still satisfy the assumption({\bf A}),
hence we have $\lambda_{n,\varepsilon}(\bar{a},\bar{b}) \leq
\lambda_{n}(\bar{a},\bar{b})+\frac{\varepsilon}{p\delta}$, Since
$(\lambda_{n,\varepsilon}(\bar{a},\bar{b}))_{n}\nearrow \infty$,
$\displaystyle\lim_{n \rightarrow
\infty}\lambda_{n}(\bar{a},\bar{b})=+\infty$.

We claim that $\lambda_{n}(a,b) \geq \lambda_{n}(\bar{a},\bar{b})$,
so we get $\displaystyle\lim_{n\rightarrow
\infty}\lambda_{n}=\displaystyle\lim_{n\rightarrow
\infty}\lambda_{n}(a,b)=+\infty$.

Suppose that $K$ is a compact  symmetric set such that $i(K) \geq n$
and $F(\bar{a},\bar{b})(u) = 1$, $\forall u \in K$. Then the map
$\Psi: K \rightarrow \Psi(K)$ defined by $u\mapsto
\frac{u}{F(a,b)(u)}$, is an odd homeomorphism, and $F(a,b)(w) = 1$,
$\forall w\in \Psi(K)$. Since $\bar{a} \geq a$ and $\bar{b} \geq b$,
we have $F(a,b)(u) \leq 1$, $\forall u\in K$. So
$\displaystyle\sup_{u\in \Psi(K)} G(u) \geq \displaystyle\sup_{w\in
K}G(w)$ and $i(\Psi(K)) = i(K) \geq n$. Therefore we have
$$\displaystyle\sup_{u\in \Psi(K)}G(u) \geq
\lambda_{n}(\bar{a},\bar{b}).$$ But any set in
$\mathcal{F}^{c}_{n}(a,b)$ can be write as the image of a set in
$\mathcal{F}^{c}_{n}(\bar{a},\bar{b})$ under the map $\Psi$, so we
get $\lambda_{n}(a,b) \geq \lambda_{n}(\bar{a},\bar{b})$. \qed
\begin{lem}\label{l3.10}
$(\lambda_{n})_{n}$ is sequence of eigenvalues associated to the
problem $G'(u) = \lambda F'(u)$.
\end{lem}
\pf: Fix $n\in \mathbf{N}^{*}$, let $\varepsilon = \frac{1}{k}$,
$k\in \mathbf{N}^{*}$. From the above discussion there exists a
sequence $(u_{k})_{k\in \mathbf{N}^{*}} $ of eigenfunctions
associated to $(\lambda_{n,\frac{1}{k}})_{k}$ satisfying
$G(u_{k})+\|u_{k}\|^{p}_{p}=1$. Hence $(u_{k})_{k}$ is bounded in
$X$, thus, up to a subsequence, $(u_{k})_{k}$ converges weakly in
$X$ to some $u\in X$. Since $u_{k}$ satisfies
$G'(u_{k})+\frac{1}{k}|u_{k}|^{p-2}u_{k}=
\lambda_{n,\frac{1}{k}}F'(u_{k})$, from Lemma \ref{l3.3}, we have
$|u_{k}|^{p-2}u_{k} \rightarrow |u|^{p-2}u$, $F'(u_{k}) \rightarrow
F'(u)$ in $X$ as $k \rightarrow \infty$, so $G'(u_{k}) \rightarrow
\lambda_{n}F'(u)$ as $k \rightarrow \infty$, $G'(u_{k}) +
|u_{k}|^{p-2}u_{k} \rightarrow \lambda_{n}F'(u) + |u|^{p-2}u$ in $X$
as $k \rightarrow \infty$. Thus there holds $$\langle G'(u_{k}) +
|u_{k}|^{p-2}u_{k}, u_{k} - u \rangle \rightarrow 0.$$ By Lemma \ref{l3.2}
with $\varepsilon = 1$, we have $u_{k} \rightarrow u$, so
$G'(u)=\lambda_{n}F'(u)$ and $(\lambda_{n}, u)$ is a solution of the
problem $G'(u) = \lambda F'(u)$.\qed
\subsection{Existence results}
The following theorems is direct consequence of subsection $3.1$.

\begin{thm}\label{l3.11}
(Existence of eigenvalue sequence for $S(\Omega)_{\varepsilon}$).
Let $F$ and $G_{\varepsilon}$ be defined in section $3.1$ with $a
\equiv 0$, $b \equiv V$ and $\beta(x) \equiv 0 $. Let $X$ be
$W^{1,p}(\Omega)$, then there exist a nondecreasing sequence of
nonnegative eigenvalues $\{\lambda_{n,\varepsilon}\}$ of
(\ref{e3.1})(when $\varepsilon = 0$, set $\lambda_{n,\varepsilon} =
\lambda_{n}$), that is , the eigenvalues of
$S(\Omega)_{\varepsilon}$, moreover, this sequence is divergent.
\end{thm}

\begin{thm}\label{l3.12}
(Existence of eigenvalue sequence for $P(\Omega)_{\varepsilon}$) Let
$F$ and $G_{\varepsilon}$ be defined in section $3.1$ with $a \equiv
V$, $b \equiv 0$ and $\beta \equiv 0$. Let $X$ be
$W^{1,p}_{0}(\Omega)\oplus \mathbf{R}$, then there exist a
nondecreasing sequence of nonnegative eigenvalues $\{\lambda_{n,
\varepsilon}\}$ of (\ref{e3.1}) (when $\varepsilon = 0$, set
$\lambda_{n,\varepsilon} = \lambda_{n}$), that is , the eigenvalues
of $P(\Omega)_{\varepsilon}$, moreover, this sequences is divergent.
\end{thm}
\begin{thm}\label{l3.13}
(Existence of eigenvalue sequence for $N(\Omega)_{\varepsilon}$) Let
$F$ and $G_{\varepsilon}$ be defined in section $3.1$ with $a \equiv
V$, $b \equiv 0$ and $\beta \equiv 0$. Let $X$ be $W^{1,p}(\Omega)$,
then there exist a nondecreasing sequence of nonnegative eigenvalues
$\{\lambda_{n, \varepsilon}\}$ of (\ref{e3.1}) (when $\varepsilon = 0$,
set $\lambda_{n,\varepsilon} = \lambda_{n}$), that is , the
eigenvalues of $N(\Omega)_{\varepsilon}$, moreover, this sequences
is divergent.
\end{thm}
\begin{thm}\label{l3.14}
(Existence of eigenvalue sequence for $R(\Omega)_{\varepsilon}$).
Let $F$ and $G_{\varepsilon}$ be defined in section $3.1$ with $a
\equiv V$, $b \equiv 0$ and $\beta(x) \equiv \gamma(x) $. Let $X$ be
$W^{1,p}(\Omega)$, then there exist a nondecreasing sequence of
nonnegative eigenvalues $\{\lambda_{n, \varepsilon}\}$ of (\ref{e3.1})(when $\varepsilon = 0$, set $\lambda_{n,\varepsilon} =
\lambda_{n}$), that is , the eigenvalues of
$R(\Omega)_{\varepsilon}$, moreover, this sequence is divergent.
\end{thm}

\subsection{Index computation for cones}
Similar to Theorem $3.2$ in \cite{DL07}, we have:
\begin{thm}\label{l3.15}If  $\lambda_{m,\varepsilon} <
\lambda_{m+1,\varepsilon}$ for some $m\in \mathbf{N}^{*}$,then\\
$i(\{u\in X\setminus \{0\}:G_{\varepsilon}(u) \leq
\lambda_{m,\varepsilon} F(u)\})$ = $i(\{u\in X: G_{\varepsilon}(u) <
\lambda_{m+1,\varepsilon} F(u)\})$ = m.
\end{thm}
\pf: Suppose $\lambda_{m,\varepsilon} < \lambda_{m+1,\varepsilon}$.
If we set $A = \{u\in \mathcal{M}: G_{\varepsilon}(u) \leq
\lambda_{m,\varepsilon}\}$ and $B = \{u\in
\mathcal{M}:G_{\varepsilon}(u) < \lambda_{m+1, \varepsilon}\}$,
clearly, we have $i(A) \leq m$. Assume that $i(A) \leq m-1$. By (i3)
in section 2.3, there exists a symmetric neighborhood $W$ of $A$ in
$\mathcal{M}$ satisfying $i(W) = i(A)$. Notice that such a $W$ is
also a neighborhood of the critical set of
$G_{\varepsilon}|_{\mathcal{M}}$ at level $\lambda_{m,\varepsilon}$,
by the equivariant deformation theorem, there exists $\delta > 0$
and an odd continuous map $\iota:\{u\in \mathcal{M}:
G_{\varepsilon}(u) \leq \lambda_{m,\varepsilon}+\delta\} \rightarrow
\{u\in \mathcal{M}: G_{\varepsilon}(u) \leq \lambda_{m,
\varepsilon}-\delta\}\cup W = W$. It follows from (i2) in section
2.3 that $i(u\in \mathcal{M}:G_{\varepsilon}(u) \leq
\lambda_{m,\varepsilon}+\delta) \leq m-1$. This contradicts the
definition of $\lambda_{m,\varepsilon}$ and the monotonicity of the
cohomological index. By the invariance of the cohomological index
under odd homeomorphism, we have $$i(\{u\in X\setminus
\{0\}:G_{\varepsilon}(u) \leq \lambda_{m,\varepsilon} F(u)\}) = m.$$

By the monotonicity of the cohomological index, we have $i(B) \geq
m$. Assume that $i(B) \geq m+1$. From the proof of Lemma \ref{l3.8},
there exists a symmetric, compact subset $K$ of $B$ with $i(K) \geq
m+1$. Since $\max\{G_{\varepsilon}(u): u\in K\} <
\lambda_{m+1,\varepsilon}$, this contradicts to Lemma \ref{l3.8}. By the
invariance of the cohomological index under odd homeomorphism, we
have $i(\{u\in X: G_{\varepsilon}(u) < \lambda_{m+1,\varepsilon}
F(u)\}) = m $.\qed
\section{Proof of the main theorem}
In this section, we assume that $\varepsilon > 0$.
\subsection{Proof of Theorem 1.1 and Theorem 1.1$'$}
We consider the $C^{1}$ functional $f_{\varepsilon}: X =
W^{1,p}(\Omega) \rightarrow \mathbf{R}$ defined by
\begin{eqnarray*}
f_{\varepsilon}(u) = \frac{1}{p}\int_{\Omega}(|\nabla u|^{p} +
\varepsilon|u|^{p}){\rm d}x -
\frac{\lambda}{p}\int_{\partial\Omega}V|u|^{p}{\rm d}S_x -
\int_{\partial\Omega}H(x,u){\rm d}S_x.
\end{eqnarray*}
It is clear that critical points of $f_{\varepsilon}$ are weak
solutions of (\ref{e1.1}).

In this case, \begin{eqnarray*} G_{\varepsilon}(u) =
\frac{1}{p}\int_{\Omega}(|\nabla u|^{p}+\varepsilon|u|^{p}){\rm d}x
\end{eqnarray*}
\begin{eqnarray*}
F(u) = \frac{1}{p}\int_{\partial\Omega}V(x)|u|^{p}{\rm d}S_x,
\end{eqnarray*}
Hence
\begin{eqnarray*}
f_{\varepsilon}(u) = G_{\varepsilon}(u)-\lambda F(u)-
\int_{\partial\Omega}H(x,u){\rm d}S_x.
\end{eqnarray*}

We follow the line of \cite{DL07}.
\begin{lem}\label{l4.1}
By (h1) and (h2), we have, $ \frac{\int_{\partial\Omega}H(x,u){\rm
d}S_x}{\|u\|_{\varepsilon }^{p}} \rightarrow 0$ as
$\|u\|_{\varepsilon} \rightarrow 0$.
\end{lem}
\pf: {\bf Case 1:} $p < N$. Set
\begin{eqnarray*}
H_{0}(x,s)= \left \{
\begin{array}
{ll} \frac{H(x,s)}{|s|^{p}},
&\textrm{  if $s \neq 0$},\\
0,
&\textrm{  if $s = 0$},\\
\end{array}
\right.
\end{eqnarray*}
from (h1) and (h2) we have that $H_{0}$ is a Carath\'{e}odory
function satisfying
\begin{eqnarray*}
|H_{0}(x,s)| \leq \frac{1}{p}a_{\varepsilon}(x) +
\frac{\varepsilon}{\frac{Np-p}{N-p}}|s|^{\frac{p^{2}-p}{N-p}}.
\end{eqnarray*}
By the continuous embedding of $X$ into
$L^{\frac{Np-p}{N-p}}(\partial\Omega)$ and Lemma \ref{l2.5}, it follows
that $H_{0}(x,u)$ converges to $0$ in
$L^{\frac{N-1}{p-1}}(\partial\Omega)$ as $\|u\|_{\varepsilon}
\rightarrow 0$. Using H\"{o}lder inequality we have
\begin{eqnarray*}
\int_{\partial\Omega}|H(s,u)|{\rm d}S_x=
\int_{\partial\Omega}|H_{0}(s,u)||u|^{p}{\rm d}S_x \leq
\Big(\int_{\partial\Omega}|H_{0}(s,u)|^{\frac{N-1}{p-1}}{\rm d}S_x
\Big)^{\frac{p-1}{N-1}}
\Big(\int_{\partial\Omega}|u|^{\frac{Np-p}{N-p}}{\rm d}S_x
\Big)^{\frac{N-p}{N-1}}.
\end{eqnarray*}
Applying Sobolev embedding theorem again, the conclusion follows in
this case.

{\bf Case 2:} $p = N$. In this case, by making $q$ large enough, we
can also write $|H_{0}(x,s)| \leq \frac{1}{p}a_{\varepsilon}(x) +
\frac{\varepsilon}{q}|s|^{q-p}$ and $q-p > 1$, here $H_{0}(x,s)$ is
defined as {\bf Case 1}. By the continuous embedding of $X$ into
$L^{r(q-p)}(\partial\Omega)$, it follows from Lemma \ref{l2.5} that
$H_{0}(x,u)$ converges to $0$ in $L^{r}(\partial\Omega)$ as
$\|u\|_{\varepsilon} \rightarrow 0$. Using the H\"{o}lder inequality
we have
\begin{eqnarray*}
\int_{\partial\Omega}|H(s,u)|{\rm d}S_x  =
\int_{\partial\Omega}|H_{0}(s,u)||u|^{p}{\rm d}S_x  \leq
\Big(\int_{\partial\Omega}|H_{0}(s,u)|^{r}{\rm d}S_x
\Big)^{\frac{1}{r}}\Big(\int_{\partial\Omega}
|u|^{p\frac{r}{r-1}}{\rm d}S_x \Big)^{\frac{r-1}{r}}.
\end{eqnarray*} Applying Sobolev embedding theorem again, the conclusion follows in this case.

{\bf Case 3:} $p > N$. In this case, we can also write $|H_{0}(x,s)|
\leq \frac{1}{p}a_{S}(x) + \frac{\varepsilon}{q}|s|^{q-p}$, for $|s|
\leq S$ and $q-p
> 1$, here $H_{0}(x,s)$ is defined as {\bf Case 1}. By Sobolev embedding theorem,
we can also assume that $\|u\|_{C^{0}(\partial\Omega)} < S$ for some
$S > 0$ when $\|u\|_{\varepsilon}$ is small. Since $X$ continuously
embeds into $L^{(q-p)}(\partial\Omega)$ and from Lemma \ref{l2.5}, Remark
\ref{l2.6}, we can deduce that $H_{0}(x,u)$ goes to $0$ in
$L^{1}(\partial\Omega)$ as $\|u\|_{\varepsilon} \rightarrow 0 $.
Using the H\"{o}lder inequality we have
\begin{eqnarray*}
\int_{\partial\Omega}|H(s,u)|{\rm d}S_x  =
\int_{\partial\Omega}|H_{0}(s,u)||u|^{p}{\rm d}S_x  \leq
(\int_{\partial\Omega}|H_{0}(s,u)|{\rm d}S_x
)\|u\|_{L^{\infty}(\partial\Omega)}^{p}.
\end{eqnarray*}
Applying Sobolev embedding theorem again, the conclusion follows in
this case.\qed
\begin{lem}\label{l4.2}
If there exists $b > 0$ and $(u_{k})$ in $X$ such that
$\|u_{k}\|_{\varepsilon} \rightarrow \infty$ and
$\int_{\Omega}(|\nabla u_{k}|^{p}+\varepsilon|u_{k}|^{p}){\rm d}x
\leq b\int_{\partial\Omega}V(x)|u_{k}|^{p}{\rm d}S_x$. Then from
(h2) and (h4) we have $\frac{\int_{\partial\Omega}H(x,u_{k}){\rm
d}S_x}{\|u_{k}\|_{\varepsilon}^{p}} \rightarrow +\infty$.
\end{lem}
\pf: Set $v_{k}=\frac{u_{k}}{\|u_{k}\|_{\varepsilon}}$, then, up to
a subsequence, $(v_{k})$ converges to some $v$ weakly in $X$ and
a.e.in $\partial\Omega$. By Lemma \ref{l3.4}, it follows that
$b\int_{\partial\Omega}V|v|^{p}ds \geq 1$. So $|v| \neq 0$ on a set
with positive measure. Thus from (h2) we have
\begin{eqnarray*}
\lim_{k \rightarrow
\infty}\frac{H(s,u_{k}(s))}{\|u_{k}\|_{\varepsilon }^{p}} =
\lim_{k\rightarrow\infty}\frac{H(s,\|u_{k}\|_{\varepsilon}v_{k}(s))}
{\|u_{k}\|_{\varepsilon}^{p}|v_{k}(s)|^{p}} |v_{k}(s)|^{p}=+\infty
\end{eqnarray*}
on a set with positive measure. By (h4) we can apply Fatou's lemma
to the sequence
$(\frac{H(s,u_{k})}{\|u_{k}\|_{\varepsilon}^{p}})_{k}$ and the
assertion follows.\qed
\begin{lem}\label{l4.3}
Suppose (h1) is satisfied. The map $T: X \rightarrow X^{*}$ defined
by $T(u)(v) = \int_{\partial\Omega}h(x,u)v{\rm d}S_x$ is
weak-to-strong continuous.
\end{lem}
\pf: If $p < N$, we set $\alpha =
\frac{\frac{Np-p}{N-p}}{\frac{Np-p}{N-p}-1} = \frac{Np-p}{Np-N}$.
Let $(u_{k})$ be a sequence weakly convergent to $u$ in $X$, then
$(u_{k})$ is bounded in $L^{\frac{Np-p}{N-p}}(\partial\Omega)$ and
up to subsequence, converges to $u$ a.e.in $\partial\Omega$. By (h1)
and Young's inequality we have
\begin{eqnarray*}|h(x,s)| &\leq& a_{\varepsilon}(x)|s|^{p-1} +
\varepsilon |s|^{\frac{Np-p}{N-p}-1} \\
&\leq & \frac{\alpha
(p-1)}{N-1}(\frac{a_{\varepsilon}}{\varepsilon})^{\frac{N-1}{\alpha
(p-1)}} +
\frac{p-1}{\frac{Np-p}{N-p}-1}\varepsilon^{\frac{\frac{Np-p}{N-p}-1}{p-1}}|s|^{\frac{Np-p}{N-p}-1}
+ \varepsilon |s|^{\frac{Np-p}{N-p}-1}.\end{eqnarray*} From Lemma
\ref{l2.5}, $(h(x,u_{k}))$ is convergent to $h(x,u)$ strongly in
$L^{\alpha}(\partial\Omega)$, hence strongly in $X^{*}$.

If $p = N$, in this case, by making $q$ large enough, we may assume
that for every $\varepsilon
> 0$, there exists $a_{\varepsilon}\in L^{r}(\partial\Omega)$ such
that $|h(x,s)| \leq a_{\varepsilon} (x)|s|^{p-1} + \varepsilon
|s|^{q - 1}$ and $q -p > 1$, $\frac{q-p}{q-1}r > 1$. Let $(u_{k})$
be a sequence weakly convergent to $u$ in $X$. Then $(u_{k})$ is
bounded in $L^{r(q - p)}(\partial\Omega)$ and up to subsequence,
converges to $u$ a.e.in $\partial\Omega$. By (h1) and Young's
inequality we have
\begin{eqnarray*}|h(x,s)| \leq a_{\varepsilon}(x)|s|^{p-1} +
\varepsilon |s|^{q-1} \leq \frac{
(q-p)}{q-1}(\frac{a_{\varepsilon}}{\varepsilon})^{\frac{q - 1}{ q-
p}} + \frac{p-1}{q-1}\varepsilon^{\frac{q-1}{p-1}}|s|^{q-1} +
\varepsilon |s|^{q-1}.\end{eqnarray*} From Lemma \ref{l2.5}, we have
$(h(x,u_{k}))$ is convergent to $h(x,u)$ strongly in
$L^{\frac{q-p}{q-1}r}(\partial\Omega)$, hence strongly in $X^{*}$.

If $p > N$, let $(u_{k})$ be a sequence weakly convergent to $u$ in
$X$. Then by Sobolev embedding thoerem, $(u_{k})$ converges to $u$
uniformly in $\partial\Omega$. By (h1), we have $$|h(x,u_{k}) -
h(x,u)| \leq a_{S}(x)(|u_{k}|^{p-1} + |u|^{p-1})$$ for some $S > 0$.
Applying Fatou's lemma to the sequence $a_{S}(x)(|u_{k}|^{p-1} +
|u|^{p-1}) - |h(x,u_{k}) - h(x,u)|$, we obtain
\begin{eqnarray*}2\int_{\partial\Omega}a_{S}(s)|u|^{p-1}{\rm d}S_x &\leq&
\liminf_{k \rightarrow
\infty}\int_{\partial\Omega}[a_{S}(s)(|u_{k}|^{p-1} + |u|^{p-1}) -
|h(s,u_{k}) - h(s,u)|]{\rm d}S_x\\ &\leq &
2\int_{\partial\Omega}a_{S}(s)|u|^{p-1}{\rm d}S_x - \limsup_{k
\rightarrow \infty}\int_{\partial\Omega}|h(s,u_{k}) - h(s,u)|{\rm
d}S_x.
\end{eqnarray*}
So $\limsup_{k \rightarrow \infty}\int_{\partial\Omega}|h(s,u_{k}) -
h(s,u)|{\rm d}S_x \leq 0$, that is, $h(x,u_{k})$ converges to
$h(s,u)$ in $L^{1}(\partial\Omega)$. From Sobolev embedding theorem,
we have that $h(x,u_{k})$ converges to $h(s,u)$ in $X^{*}$. \qed

\begin{lem}\label{l4.4}
Suppose (h1)--(h4) hold.  For every $\lambda \in \mathbf{R}$ and
$c\in \mathbf{R}$, the functional $f_{\varepsilon}$ satisfies
$(PS)_{c}$ condition.
\end{lem}
\pf: Let $(u_{k})_{k}$ be a sequence in $X$ satisfying
$f'_{\varepsilon}(u_{k}) \rightarrow 0$ in $X^{*}$ and
$f_{\varepsilon}(u_{k})\rightarrow c$.

{\bf Claim:} $(u_{k})$ is bounded in $X$. By contradiction, we
assume that $\|u_{k}\|_{\varepsilon} \rightarrow \infty$. From (h3)
we have
\begin{eqnarray*}&&\mu f_{\varepsilon}(u_{k}) - \langle
f'_{\varepsilon}(u_{k}),u_{k}\rangle =
(\frac{\mu}{p}-1)\int_{\Omega}(|\nabla u_{k}|^{p}+ \varepsilon
|u_{k}|^{p})dx -(\frac{\mu}{p}-1)\int_{\partial\Omega}\lambda
V|u_{k}|^{p}{\rm d}S_x\\&&+
\int_{\partial\Omega}(h(s,u_{k})u_{k}-\mu H(s,u_{k})){\rm d}S_x \geq
(\frac{\mu}{p}-1)\int_{\Omega}(|\nabla u_{k}|^{p}+ \varepsilon
|u_{k}|^{p})dx\\&& -(\frac{\mu}{p}-1)\int_{\partial\Omega}\lambda
V|u_{k}|^{p}{\rm d}S_x  - \int_{\partial\Omega}(\gamma_{0}
+\gamma_{1}|u_{k}|^{p}){\rm d}S_x.\end{eqnarray*} Since
\begin{eqnarray*}\mu f_{\varepsilon}(u_{k}) - \langle
f'_{\varepsilon}(u_{k}),u_{k}\rangle +
\int_{\partial\Omega}\gamma_{0}{\rm d}S_x \leq
\frac{1}{2}(\frac{\mu}{p}-1)\int_{\Omega}(|\nabla u_{k}|^{p} +
\varepsilon |u_{k}|^{p})dx \end{eqnarray*} for $k$ large enough,
there exists $b
> 0$ such that $$\int_{\Omega}(|\nabla u_{k}|^{p}+ \varepsilon
|u_{k}|^{p})dx \leq \int_{\partial\Omega}(2\lambda V +
b\gamma_{1})|u_{k}|^{p}{\rm d}S_x$$ for $k$ large enough. $(2\lambda
V + b\gamma_{1})$ is still in $L^{r}(\partial\Omega)$, from Lemma
\ref{l4.2} we can deduce that
\begin{eqnarray*}
\lim_{k\rightarrow \infty}\frac{\int_{\partial\Omega}H(s,u_{k}){\rm
d}S_x}{\|u_{k}\|_{\varepsilon}^{p}} = +\infty.\end{eqnarray*}
Moreover, by Sobolev embedding theorem, we have
$\int_{\partial\Omega}V|u|^{p}{\rm d}S_x \leq
C\|V\|_{r}\|u\|_{\varepsilon}^{p}$. Therefore
\begin{eqnarray*} 0 =
\lim_{k\rightarrow
\infty}\frac{f_{\varepsilon}(u_{k})}{\|u_{k}\|_{\varepsilon}} =
\frac{1}{p} -
\lim_{k\rightarrow\infty}(\frac{\lambda\int_{\partial\Omega}V|u_{k}|^{p}{\rm
d}S_x}{p\|u_{k}\|_{\varepsilon}^{p}} + \frac{\int_{\partial\Omega}
H(s,u_{k}){\rm d}S_x}{\|u_{k}\|_{\varepsilon}^{p}}) = -\infty.
\end{eqnarray*}
It is a contradiction, so $(u_{k})$ is bounded in $X$.

Actually, $f'_{\varepsilon}(u_{k}) = G'_{\varepsilon}(u_{k}) -
\lambda F'(u_{k}) - T(u_{k})$, here $T : X \rightarrow X^{*}$ is
defined in Lemma \ref{l4.3}. By Lemma \ref{l3.3} and \ref{l4.3}, we have, up to subsequence,
$F'(u_{k})$ and $T(u_{k})$ converge, so $G'_{\varepsilon}(u_{k})$
converges in $X^{*}$. By Lemma \ref{l3.2}, we can deduce that $u_{k}$ has a
convergent subsequence. So we have proved the $(PS)_{c}$
condition.\qed

In order to prove the Theorem 1.1$'$, we need the following result.\\
{\bf{Lemma 4.4$'$}}  {\it Suppose (h1),(h2),(h4),(h5) hold.  For every $\lambda \in \mathbf{R}$,
$f_{\varepsilon}$ satisfies
the Cerami condition.
}\\
\pf: Let $(u_{k})_{k}$ be a sequence in $X$ satisfying
$(1+\|u_{k}\|_{\varepsilon})f'_{\varepsilon}(u_{k}) \rightarrow 0$ in $X^{*}$ and
$f_{\varepsilon}(u_{k})\rightarrow c$.

{\bf Claim:} $(u_{k})$ is bounded in $X$. Otherwise, if $\|u_k\|_{\varepsilon}\to \infty$,  we consider $w_{k}:=\frac{u_{k}}{\|u_{k}\|_{\varepsilon}}$. Then,  up to subsequence,  we get $w_k \rightharpoonup w$ in $X$ and $w_k(x)\rightarrow w(x)$ a.e. $x\in \partial\Omega$ as $k\to \infty$.
If $w \neq 0$ in $X$,  since $f_{\varepsilon}'(u_{k})u_{k} \rightarrow 0$, that is to say
\begin{equation}\label{e4.1}
\int_{\Omega}(|\nabla u_{k}|^{p} + \varepsilon|u_{k}|^{p}){\rm d}x - \lambda\int_{\partial\Omega}V(x)|u_{k}|^{p}{\rm d}S_x -
\int_{\partial\Omega}h(x,u_{k})u_{k}{\rm d}S_x \rightarrow 0,
\end{equation}
by Schwartz inequality and Sobolev embedding theorem, we have
\begin{eqnarray*}
\frac{|\int_{\partial\Omega}V(x)|u_{k}|^{p}{\rm d}S_x|}{\|u_{k}\|_{\varepsilon}^{p}}\leq C\|V\|_{r},
\end{eqnarray*}
so by dividing the left hand side of (\ref{e4.1}) with $\|u_k\|_{\varepsilon}^p$ there holds \begin{equation}\label{e4.2}\left|\int_{\partial\Omega}\frac{h(x, u_{k})u_{k}}{\|u_{k}\|_{\varepsilon}^{p}}{\rm d}S_x\right| \leq C.\end{equation} On the other hand, by condition (h5), we have $h(x,s)s \geq H(x,s)$, so by condition (h2), $\displaystyle\lim_{|s| \rightarrow
\infty}\frac{h(x,s)s}{|s|^{p}} = +\infty$.
By Fatou's lemma, we have
\begin{eqnarray*}
\int_{\partial\Omega}\frac{h(x, u_{k})u_{k}}{\|u_{k}\|_{\varepsilon}^{p}}{\rm d}S_x = \int_{\{w_{k} \neq 0\}}|w_{k}|^{p}\frac{h(x, u_{k})u_{k}}{|u_{k}|^{p}}{\rm d}S_x \rightarrow \infty,
\end{eqnarray*}
this contradicts to (\ref{e4.2}).

If $w=0$ in $X$,  inspired by \cite{J99}, we choose $t_k\in [0,1]$ such that $f_{\varepsilon}(t_{k}u_{k}):= \displaystyle\max_{t\in [0, 1]}f_{\varepsilon}(tu_{k})$. For any $\beta> 0$ and $\tilde{w}_{k}:=(2p\beta)^{1/p}w_{k}$,    by Lemma 3.3 and Lemma 4.3 we have that
\begin{eqnarray*}
f_{\varepsilon}(t_{k}u_{k}) \geq f_{\varepsilon}(\tilde{w}_{k}) = 2\beta - \frac{\lambda}{p} \int_{\partial\Omega}V(x)|\tilde{w}_{k}|^{p}{\rm d}S_x - \int_{\partial\Omega}H(x, \tilde{w}_{k}){\rm d}S_x  \geq \beta,
\end{eqnarray*}
when $k$ is large enough, this implies that
\begin{equation}\label{e4.3}\lim\limits_{k \rightarrow \infty}f_{\varepsilon}(t_{k}u_{k}) = \infty.\end{equation} Since
$f_{\varepsilon}(0)=0,\;f_{\varepsilon}(u_k)\to c$,  we have $t_{k}\in (0, 1)$.
By the definition of $t_k$,
\begin{equation}\label{e4.4}\langle f_{\varepsilon}'(t_{k}u_{k}), t_{k}u_{k}\rangle = 0.\end{equation}
From (\ref{e4.3}), (\ref{e4.4}), we have
\begin{eqnarray*}
f_{\varepsilon}(t_{k}u_{k})-\frac 1p\langle f_{\varepsilon}'(t_{k}u_{k}), t_{k}u_{k}\rangle =\int_{\partial\Omega}\left(\frac{1}{p}h(x, t_{k}u_{k})t_{k}u_{k} - H(x, t_{k}u_{k})\right){\rm d}S_x \rightarrow \infty.
\end{eqnarray*}
By  (h3),  there exists $\theta\ge 1$ such that
\begin{equation}\label{e4.5}
\int_{\partial\Omega}(\frac{1}{p}h(x, u_{k})u_{k} - H(x, u_{k})){\rm d}S_x
\geq \frac{1}{\theta}\int_{\partial\Omega}(\frac{1}{p}h(x, t_{k}u_{k})t_{k}u_{k} - H(x, t_{k}u_{k})){\rm d}S_x \rightarrow \infty.
\end{equation}
On the other hand,
\begin{equation}\label{e4.6}
\int_{\partial\Omega}(\frac{1}{p}h(x, u_{k})u_{k} - H(x, u_{k})){\rm d}S_x
= f_{\varepsilon}(u_{k}) - \frac{1}{p}\langle f_{\varepsilon}'(u_{k}),u_{k}\rangle\to c.
\end{equation}
(\ref{e4.5}) and (\ref{e4.6}) are contradiction. Hence $\{u_{k}\}$ is bounded in $X$. So up to a subsequence, we can assume that $u_{k}\rightharpoonup u$ for some $X$.

The same reason as Lemma \ref{l4.4}, we can prove that $\{u_{k}\}$ have a convergent subsequence. So $f_{\varepsilon}$ satisfies the Cerami condition. \qed \\
{\bf Proof of Theorem 1.1}: Replacing $(\lambda, V)$ with
$(-\lambda, -V)$, we can assume that $\lambda \geq 0$.

{\bf Case 1}: $meas\{x\in \partial\Omega: V(x) > 0\}> 0$ (by Theorem
\ref{l3.11}, $S(\Omega)_{\varepsilon}$ has a divergent sequence
$(\lambda_{m,\varepsilon})_{m}$ of eigenvalues), $\lambda \geq
\lambda_{1,\varepsilon}$.

Since the sequence $(\lambda_{m,\varepsilon})_{m}$ is divergent,
there exist $m \geq 1$ such that $\lambda_{m,\varepsilon} \leq
\lambda <\lambda_{m+1,\varepsilon}$. Define
\begin{eqnarray*}
C_{-} = \{u\in X: G_{\varepsilon}(u) \leq \lambda_{m,\varepsilon}
F(u)\},\\
 C_{+} = \{u\in X: G_{\varepsilon}(u) \geq
\lambda_{m+1,\varepsilon} F(u)\},
\end{eqnarray*}
we have that $C_{-}$, $C_{+}$ are two symmetric closed cones in $X$
with $C_{-} \cap C_{+} = \{0\}$.

By Theorem \ref{l3.15} we have that $i(C_{-}\setminus\{0\}) = i(X\setminus
C_{+}) = m$.

Since $\lambda < \lambda_{m+1,\varepsilon}$, by Lemma \ref{l4.1} there
exist $r_{+}
> 0$ and $\alpha > 0$ such that $f_{\varepsilon}(u)
> \alpha$ for $u\in C_{+}$ and $\|u\|_{\varepsilon} = r_{+}$.
Since $\lambda \geq \lambda_{m,\varepsilon}$, by (h4) we have
$f_{\varepsilon}(u) \leq 0$ for every $u\in C_{-}$.

Let $e\in X\setminus C_{-}$, we define another norm on $X$ by
$\|u\|_{V}:=(\int_{\partial\Omega}(|V|+1)|u|^{p}{\rm d}S_x)^{1/p}$.
If $u\in C_{-}$ and $t > 0$, then \begin{eqnarray*}\|u +
te\|_{\varepsilon} = t\|\frac{u}{t} + e\|_{\varepsilon} \leq
t(\|\frac{u}{t}\|_{\varepsilon} + \|e\|_{\varepsilon}) \leq
t(C\|\frac{u}{t}\|_{V} +
\frac{\|e\|_{\varepsilon}}{\|e\|_{V}}\|e\|_{V}) \leq C
t(\|\frac{u}{t}\|_{V} + \|e\|_{V}).\end{eqnarray*} Notice that
$C_{-}$ is also closed in $X$ with respect to the norm
$\|\cdot\|_{V}$, by Proposition 2.12 in \cite{DL07}, there exists
$\beta \geq 1$ such that $\|\frac{u}{t}\|_{V} + \|e\|_{V} \leq \beta
\|\frac{u}{t} + e\|_{V}$. Hence, $\|u + te\|_{\varepsilon} \leq b
\|u + te\|_{V}$ for every $u\in C_{-}$, $t \geq 0$ and some $b > 0$.
Thus from Lemma \ref{l4.2} we have that
$\frac{\int_{\partial\Omega}H(s,u_{k}){\rm
d}S_x}{\|u_{k}\|_{\varepsilon}^{p}}\rightarrow +\infty$ for
$\|u_{k}\|_{\varepsilon}\rightarrow +\infty$ and $u_{k}\in
C_{-}+\mathbf{R}^{+}e$. So there exists $r_{-}
> r_{+}$ such that $f_{\varepsilon}(u) \leq 0$ for $u\in C_{-} +
\mathbf{R}^{+}e$ and $\|u\|_{\varepsilon} \geq r_{-}$.

If we define $D_{-}$, $S_{+}$, $Q$, $H$ as Lemma \ref{l2.3}, then
$f_{\varepsilon}$ is bounded on $Q$, $f_{\varepsilon}(u) \leq 0$ for
every $u\in D_{-}\cup H$ and $f_{\varepsilon}(u) \geq \alpha
> 0$ for every $u\in S_{+}$. With Lemma \ref{l4.4}, it follows
that $f_{\varepsilon}$ has a critical value $c \geq \alpha
> 0$. Hence $u$ is a nontrivial weak solution of (\ref{e1.1}).

{\bf Case 2:} $meas\{x\in \partial\Omega: V(x) > 0\}> 0$, $0 \leq
\lambda < \lambda_{1,\varepsilon}$ or $meas\{x\in \partial\Omega:
V(x) > 0\} = 0$, $\lambda \geq 0$. We set $C_{-} = \{0\}$, $C_{+} =
X$ and the proof is similar.\qed\\
{\bf Proof of theorem 1.1$'$}:
The process is the same as the proof of Theorem 1.1. With the aid of the remark after Lemma 2.3, we use Lemma 4.4$'$ instead of Lemma \ref{l4.4}.\qed
\subsection{Proof of Theorem 1.2 and Theorem 1.2$'$}
For problem (\ref{e1.2}), we consider the $C^{1}$ functional
$f_{\varepsilon}: X = W^{1,p}_{0}(\Omega)\oplus\mathbf{R}
\rightarrow \mathbf{R}$ defined by
\begin{eqnarray*}
f_{\varepsilon}(u) &=& \frac{1}{p}\int_{\Omega}(|\nabla u|^{p} +
\varepsilon|u|^{p}){\rm d}x -
\frac{\lambda}{p}\int_{\Omega}V(x)|u|^{p}{\rm d}x -
\int_{\Omega}H(x,u){\rm d}x\\
&=& G_{\varepsilon}(u)-\lambda F(u)-\int_{\Omega}H(x,u){\rm d}x
~~(see~ Theorem ~\ref{l3.11}).
\end{eqnarray*}
For problem (\ref{e1.3}), we consider the $C^{1}$ functional
$f_{\varepsilon}: X = W^{1,p}(\Omega) \rightarrow \mathbf{R}$
defined by
\begin{eqnarray*}
f_{\varepsilon}(u) &=& \frac{1}{p}\int_{\Omega}(|\nabla u|^{p} +
\varepsilon|u|^{p}){\rm d}x -
\frac{\lambda}{p}\int_{\Omega}V|u|^{p}{\rm d}x -
\int_{\Omega}H(x,u){\rm d}x\\
&=& G_{\varepsilon}(u)-\lambda F(u)-\int_{\Omega}H(x,u){\rm d}x
~~(see~ Theorem ~\ref{l3.12}).
\end{eqnarray*}
For problem (\ref{e1.4}), we consider the $C^{1}$ functional
$f_{\varepsilon}: X = W^{1,p}(\Omega) \rightarrow \mathbf{R}$
defined by
\begin{eqnarray*}
f_{\varepsilon}(u) &=& \frac{1}{p}\int_{\Omega}(|\nabla u|^{p} +
\varepsilon|u|^{p}){\rm d}x
+\frac{1}{p}\int_{\partial\Omega}\beta|u|^{p}{\rm d}S_x-
\frac{\lambda}{p}\int_{\Omega}V|u|^{p}{\rm d}x - \int_{\Omega}H(x,u){\rm d}x\\
&=& G_{\varepsilon}(u)-\lambda F(u)-\int_{\Omega}H(x,u){\rm d}x
~~(see~ Theorem ~\ref{l3.13}).
\end{eqnarray*}
It is clear that critical points of $f_{\varepsilon}$ are weak
solutions of (\ref{e1.2}), (\ref{e1.3}), (\ref{e1.4}), respectively. The following
lemmas are needed in the proofs of Theorem \ref{t1.2} and 1.2$'$. Their proofs are similar to
the proofs of Lemma \ref{l4.1}, \ref{l4.2}, \ref{l4.3}, \ref{l4.4}, see also \cite{DL07}.
In the following Lemmas 4.5-4.8, we always assume that (h1)-(h4) with $\partial\Omega$
replaced by $\Omega$ and $p^{*}=\frac{Np}{N-p}$.
\begin{lem}\label{l4.5}
There holds $ \frac{\int_{\Omega}H(x,u){\rm d}x}{\|u\|_{\varepsilon
}^{p}} \rightarrow 0$ as $\|u\|_{\varepsilon} \rightarrow 0$.
\end{lem}
\begin{lem}\label{l4.6}
If there exists $b > 0$ and $(u_{k})$ in $X$ such that
$\|u_{k}\|_{\varepsilon} \rightarrow \infty$ and
$\int_{\Omega}(|\nabla u_{k}|^{p}+\varepsilon|u_{k}|^{p}){\rm d}x
\leq b\int_{\Omega}V|u_{k}|^{p}{\rm d}x$. Then from we have
$\frac{\int_{\Omega}H(x,u_{k}){\rm d}x}{\|u_{k}\|_{\varepsilon}^{p}}
\rightarrow +\infty$.
\end{lem}
\begin{lem}\label{l4.7}
The map $T: X \rightarrow X^{*}$ defined by $T(u)(v) =
\int_{\Omega}h(x,u)v{\rm d}x$ is weak-to-strong continuous.
\end{lem}
\begin{lem}\label{l4.8}
For every $\lambda \in \mathbf{R}$ and $c\in \mathbf{R}$, the
functional $f_{\varepsilon}$ satisfies $(PS)_{c}$.
\end{lem}

Under the conditions of Theorem 1.2$'$, we have the following result.\\
{\bf{Lemma 4.8$'$}} {\it  For every $\lambda \in \mathbf{R}$, $f_{\varepsilon}$ satisfies
the Cerami condition.
}\\
{\bf Proof of Theorem 1.2}: The proof is  similar to the proof of
Theorem $1.1$ by using Lemmas 4.5-4.8. We omit the details here.\qed\\
{\bf Proof of Theorem 1.2$'$}: The proof is  similar to the proof of
Theorem 1.1 by using Lemmas 4.5-4.7 and 4.8$'$. We omit the details here.\qed

{\bf \noindent \large Acknowledgments} The authors thank  Siham El
Habib for useful discussions.

\end{document}